\documentclass[11pt]{article}
\usepackage{amsthm, amsmath, amssymb, amsfonts, url, booktabs, tikz, setspace, fancyhdr, amsbsy}
\usepackage[margin = 0.8in]{geometry}
\usepackage{hyperref}       % hyperlinks
\usepackage{url}            % simple URL typesetting
\usepackage{booktabs}       % professional-quality tables
\usepackage{nicefrac}       % compact symbols for 1/2, etc.
\usepackage{microtype}      % microtypography
\usepackage{xcolor}         % colors
\usepackage{amsfonts,amssymb,bm}
\usepackage{graphicx}
\usepackage{multirow}
\usepackage{wasysym}
\usepackage{wrapfig}
\usepackage{epstopdf}
\usepackage{ulem}
\usepackage{algorithm}
\usepackage{algpseudocode}
\ifpdf
  \DeclareGraphicsExtensions{.eps,.pdf,.png,.jpg}
\else
  \DeclareGraphicsExtensions{.eps}
\fi
\newcommand{\RNum}[1]{\uppercase\expandafter{\romannumeral #1\relax}}

\newcommand{\R}{\mathbb{R}}

\newcommand{\xx}{\boldsymbol{x}}
\newcommand{\aaa}{\boldsymbol{a}}
\newcommand{\dx}{\,\mathrm{d}\boldsymbol{x}}

\newcommand{\dom}{\,\mathrm{d}{\boldsymbol{\omega}}}

\newcommand{\NN}{\mathcal{N}}
\newcommand{\LL}{\mathcal{L}}
\newcommand{\om}{\boldsymbol{\omega}}

\newtheorem{theorem}{Theorem}[section]

\newtheorem{lemma}[theorem]{Lemma}

\theoremstyle{definition}
\newtheorem{definition}[theorem]{Definition}

\newtheorem{remark}[theorem]{Remark}

\DeclareMathOperator*{\argmin}{arg\,min}

\numberwithin{equation}{section}

\def\ocirc#1{\ifmmode\setbox0=\hbox{$#1$}\dimen0=\ht0
    \advance\dimen0 by1pt\rlap{\hbox to\wd0{\hss\raise\dimen0
    \hbox{\hskip.2em$\scriptscriptstyle\circ$}\hss}}#1\else
    {\accent"17 #1}\fi}

\begin{document}

\title{Finite Element Operator Network for Solving Elliptic-type\\ Parametric PDEs}

\author{Youngjoon Hong\thanks{Department of Mathematical Sciences, Seould National University, Seoul, Republic of Korea.\\Email: \tt{hongyj@kaist.ac.kr}} , ~Seungchan Ko\thanks{Department of Mathematics, Inha University, Incheon, Republic of Korea. \\Email: \tt{scko@inha.ac.kr}}
,~and ~Jaeyong Lee\thanks{Department of AI, Chung-Ang University, Seoul, Republic of Korea. \\Email: \tt{jaeyong@cau.ac.kr}}}

\date{~}

\maketitle

~\vspace{-1.5cm}

\begin{abstract}
Partial differential equations (PDEs) underlie our understanding and prediction of natural phenomena across numerous fields, including physics, engineering, and finance. However, solving parametric PDEs is a complex task that necessitates efficient numerical methods. In this paper, we propose a novel approach for solving parametric PDEs using a Finite Element Operator Network (FEONet). Our proposed method leverages the power of deep learning in conjunction with traditional numerical methods, specifically the finite element method, to solve parametric PDEs in the absence of any paired input-output training data. We performed various experiments on several benchmark problems and confirmed that our approach has demonstrated excellent performance across various settings and environments, proving its versatility in terms of accuracy, generalization, and computational flexibility. {While our method is not meshless, the} FEONet framework shows potential for application in various fields where PDEs play a crucial role in modeling complex domains with diverse boundary conditions and singular behavior. Furthermore, we provide theoretical convergence analysis to support our approach, utilizing finite element approximation in numerical analysis.
\end{abstract}


\smallskip


\section{Introduction}
Solving partial differential equations (PDEs) is vital as they serve as the foundation for understanding and predicting the behavior of a range of natural phenomena \cite{CH-book, bender2000introduction}. From fluid dynamics, heat transfer, to electromagnetic fields, PDEs provide a mathematical framework that allows us to comprehend and model these complex systems \cite{gershenfeld1999nature,lin1988mathematics}. 
Their significance extends beyond the realm of physics and finds applications in fields such as finance, economics, and computer graphics \cite{simon1994mathematics,mortenson1999mathematics}. 
In essence, PDEs are fundamental in advancing our understanding of the world around us and have a broad impact across various industries and fields of study.

%\cite{PIDON}
%\cite{PIDON}
%%%%%%%%%%%%%%%%%%%%%%%%%%%%%%%%%%%%%%%%%%%%%%%%%%%%%%
%%%   Numerical PDEs are important and explain how they developed
%%%%%%%%%%%%%%%%%%%%%%%%%%%%%%%%%%%%%%%%%%%%%%%%%%%%%%

Numerical methods are essential for approximating solutions to PDEs when exact solutions are unattainable. This is especially significant for complex systems that defy solutions using traditional methods \cite{burden2015numerical, atkinson1991introduction}.
Various techniques, such as finite difference, finite element, and finite volume methods, are utilized to develop these numerical methods \cite{larson2013finite, strikwerda2004finite, lev02}. 
In particular, the finite element method (FEM) is widely employed in engineering and physics, where the domain is divided into small elements and the solution within each element is approximated using polynomial functions. 
The mathematical analysis of the FEM has undergone significant development, resulting in a remarkable enhancement in the method's reliability.
The FEM is particularly advantageous in handling irregular geometries and complex boundary conditions, making it a valuable tool for analyzing structures, heat transfer, fluid dynamics, and electromagnetic problems \cite{FEM1, FEM2}. 
{However, the implementation of numerical methods often comes with a significant computational cost.}

%%%%%%%%%%%%%%%%%%%%%%%%%%%%%%%%%%%%%%%%%%%%%%%%%%%%%%
%%%   ML Approach is now emerging.
%%%%%%%%%%%%%%%%%%%%%%%%%%%%%%%%%%%%%%%%%%%%%%%%%%%%%%
The use of machine learning in solving PDEs has gained significant traction in recent years. By integrating deep neural networks and statistical learning theory into numerical PDEs, a new field known as scientific machine learning has emerged, presenting novel research opportunities. The application of machine learning to PDEs can be traced back to the 1990s when neural networks were employed to approximate solutions \cite{lagaris1998artificial}. More recently, the field of physics-informed neural networks (PINNs) has been developed, where neural networks are trained to learn the underlying physics of a system, enabling us to solve PDEs and other physics-based problems \cite{PINN001,PINN007,PINN004,PINN005,PINN009}. Furthermore, this has led to active research on modified models of PINN by leveraging the advantages of it \cite{GOSWAMI2022114587, alber2019integrating,patel2021physics}. Despite their advantages, these approaches come with limitations. One prominent drawback of the PINN methods is that they are trained on a single instance of input data such as initial conditions, boundary conditions, and external force terms. As a consequence, if the input data changes, the entire training process must be repeated, posing challenges in making real-time predictions for varying input data. This limitation restricts the applicability of PINNs in dynamic and adaptive systems where the input data is subject to change.
%%%   About the OL
%To address this challenge, a new area of research called operator learning has been developed, encompassing techniques such as data-supervised or data-driven methods \cite{lu2021learning,li2021fourier,PIDON,brandstetter2022message,lienen2022learning,lu2022comprehensive}. Operator learning aims to acquire the mathematical operators governing the behavior of a physical system from data and leverage them to solve parametric PDEs. These methods rely on a pre-generated database of input-output pairs, obtained by solving the equations analytically or numerically. By utilizing this database, real-time predictions for varying input data can be made. However, generating a reliable training dataset for these methods requires a substantial number of numerical solutions obtained through extensive numerical computations, which can be computationally inefficient and time-consuming. Moreover, obtaining a sufficiently large dataset analytically or numerically can be challenging, especially for systems involving complex geometries or nonlinear equations. Another complication arises from using neural networks as the solution space, as it can introduce difficulties in imposing boundary values. Consequently, this can complicate the training process and result in less accurate solutions. Overcoming this limitation is an active area of research, with efforts focused on developing new methods to enhance the efficiency and accuracy of these data-driven approaches \cite{goswami2022physics}.
%
% Shorter version of the OL statement
%

In this regard, operator networks, an emerging field in research, address challenges by employing data-driven approaches to comprehend mathematical operators that govern physical systems. This is particularly relevant in solving parametric PDEs, as referenced in multiple studies \cite{lu2021learning,li2021fourier,PIDON,brandstetter2022message,lienen2022learning,lu2022comprehensive, miller2022neural,fanaskov2022spectral,lee2023hyperdeeponet, patel2021physics}. In particular, based upon the Universal Approximation Theorem for operators, \cite{lu2021learning} introduced the Deep Operator Network (DeepONet) architecture.
{Moreover, the authors in \cite{lu2022comprehensive} proposed an extension of the DeepONet model, called POD-DeepONet, which replaces a part of the DeepONet architecture with a prefixed basis obtained through proper orthogonal decomposition (POD) on the training data.}
These methods enable us to make fast predictions of solutions in real-time whenever the given PDE data changes, offering a new framework for solving parametric PDEs.
However, these methods rely on training data for supervised learning that consists of pre-computed (either analytically or numerically) pairs of PDE parameters (e.g., external force, boundary condition, coefficients, initial condition) and their corresponding solutions. In general, generating a reliable training dataset requires extensive numerical computations, which can be computationally inefficient and time-consuming. Obtaining a sufficiently large dataset is also challenging, especially for systems with complex geometries or nonlinear equations.
To overcome these limitations, the integration of PINNs and operator learning has led to the development of new methods, such as the Physics-Informed Neural Operator (PINO) \cite{PINO} and Physics-Informed DeepONet (PIDeepONet) \cite{PIDON}. These methods strive to leverage the strengths of both PINNs and operator learning by incorporating physical equations into the neural operator's loss function. However, they exhibit certain drawbacks, such as relatively low accuracy in handling complex geometries and difficulties in effectively addressing stiff problems like boundary layer issues, which are critical in real-world applications. Moreover, PIDeepONet tends to exhibit high generalization errors when lacking sufficient input samples.
Additionally, employing neural networks as the solution space presents challenges in imposing boundary values, potentially leading to less accurate solutions. 
%Researchers are actively engaged in developing methods to improve the efficiency and accuracy of these data-driven approaches \cite{goswami2022physics}.

% \begin{figure}[h]
% \begin{center}
% \includegraphics[width=0.8\textwidth]{domain.png}
% \end{center}
% \caption{Domain triangulation for three different complex geometries.}
% \label{fig:domain}
% \end{figure}

%%%%%%%%%%%%%%%%%%%%%%%%%%%%%%%%%%%%%%%%%%%%%%%%%%%%%%
%%%   We now introduce our method here
%%%%%%%%%%%%%%%%%%%%%%%%%%%%%%%%%%%%%%%%%%%%%%%%%%%%%%

% \marginpar{Add exp Figure ref(which was in main contribution) in here?}
To address the aforementioned limitations, this paper proposes a novel operator network for solving diverse parametric PDEs without reliance on training data. More precisely, we introduce the Finite Element Operator Network (FEONet), which utilizes finite element approximation to learn the solution operator. This approach eliminates the need for solution datasets in addressing various parametric PDEs.
\begin{figure}[t!]
\begin{center}
\includegraphics[width=\textwidth]{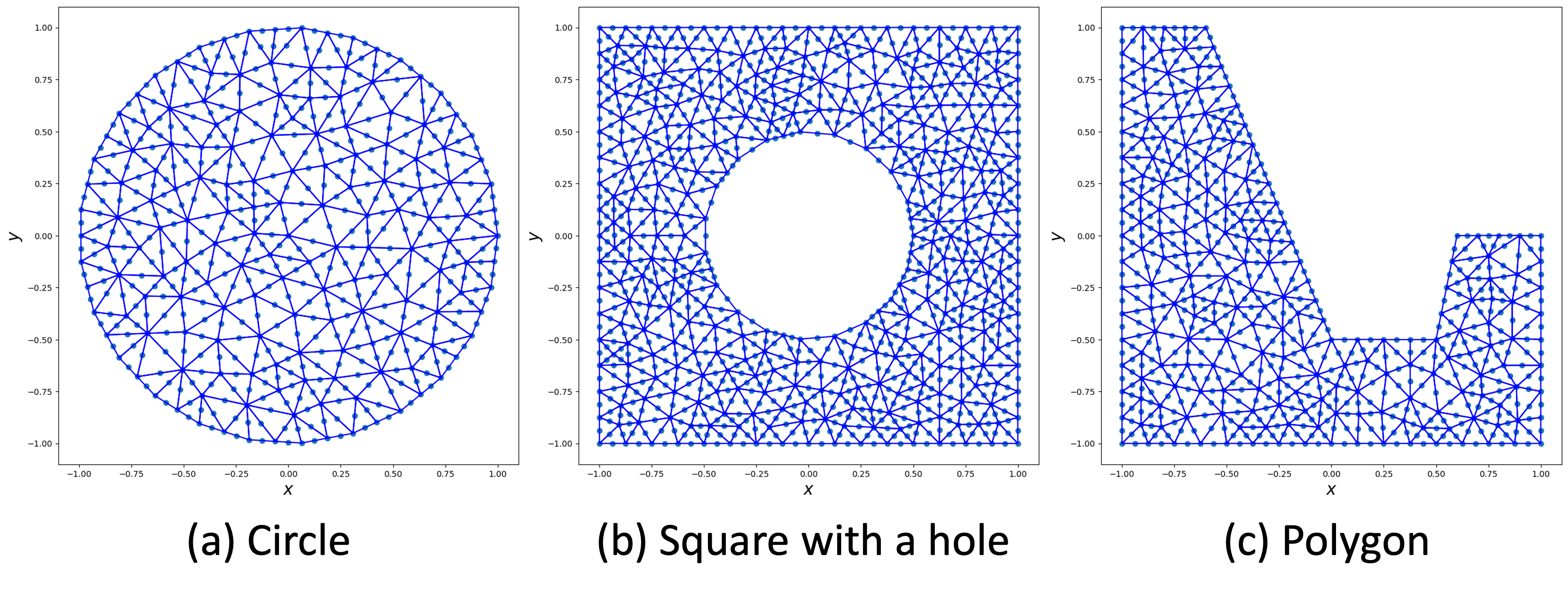}
\end{center}
\caption{Domain triangulation for three different complex geometries.}
\label{fig:domain}
\end{figure}
In the FEM framework, the numerical solution is approximated by the linear combination of nodal coefficients, $\alpha_k$, and the nodal basis, $\phi_k({\xx})$. 
The nodal basis in the FEM consists of piecewise polynomials defined by the finite set of nodes of a mesh so that
\begin{equation} \label{e:1}
{u_h}({\xx}) = \sum \alpha_k \phi_k({{\xx}}), \quad \xx \in \mathbb{R}^d,
\end{equation}
where ${u_h}$ is the FEM solution for the given PDEs.
Motivated by \eqref{e:1}, the FEONet can predict numerical solutions of the PDEs when given the initial conditions, external forcing functions, PDE coefficients, or boundary conditions as inputs. Therefore, the FEONet can learn multiple instances of the solutions of PDEs, making it a versatile approach for solving different types of PDEs in complex domains; see e.g. Figure \ref{fig:domain}.
The loss function of the FEONet is designed based on the residual quantity of the finite element approximation, motivated by the classical FEM \cite{C02, BS08}. This allows the FEONet to accurately approximate the solutions of PDEs, while also ensuring that the boundary conditions are exactly satisfied.
%To construct an approximation for the solution of PDEs, the FEONet borrows a framework from the FEM.
This approach involves inferring coefficients, denoted as $\widehat{\alpha}_k$, which are subsequently utilized to construct the linear combination $\sum\widehat{\alpha}_k\phi_k$ to approximate the solution of PDEs. 
Inheriting the capability of the FEM handling boundary conditions, the predicted solution will also satisfy the exact boundary condition. {{Note further that since we reconstruct the solution as a linear combination of the coefficients predicted by neural networks and the basis function, the FEONet scheme generates a solution as a function defined on the whole domain, not only on particular points in the domain.}} Additionally, since we have adopted the FEM as a baseline framework, we can expect that our proposed method is robust, reliable, and accurate. It is noteworthy that due to the intrinsic structure of the proposed FEONet scheme, we do not need any paired input-output training data and hence the model can be trained in an unsupervised manner. %{\color{cyan}Input sample}

%Another key contribution of this study is to introduce a novel learning architecture designed to accurately solve convection-dominated singularly perturbed problems involving strong boundary layer phenomena. These types of problems pose significant challenges for traditional machine learning and numerical methods due to sharp transitions inside thin layers caused by a small diffusive parameter. 
{Another key contribution of this study is the introduction of a novel learning architecture designed to accurately solve convection-dominated singularly perturbed problems, which are characterized by strong boundary layer phenomena. These types of problems present significant challenges for traditional machine learning and numerical methods, primarily due to the sharp transitions within thin layers caused by a small diffusive parameter. In particular, this is challenging in the realm of learning theory, where neural networks typically assume a smooth prior, making them less equipped to handle such abrupt transitions effectively.}
By adapting theory-guided methods \cite{H2020}, the network effectively captures the behavior of boundary layers, an area that has been underexplored in recent machine-learning approaches.
% One of the key contributions of this study is the development of a learning architecture that is able to accurately solve singular perturbation and boundary layer problems. 
% These types of problems, known as convection-dominated singularly perturbed problems, present significant challenges for numerical methods. 
% This is due to the fact that a small diffusive parameter produces a sharp transition inside thin layers, which can lead to large numerical errors when using traditional methods. To overcome this, semi-analytic methods such as ``enriched spaces" have been proposed and successfully applied. 
% The main component of the enriched space method is the addition of a global basis function (a so-called corrector function in the analysis) to supplement the discrete space with the corrector function. 
% Our architecture builds upon this approach by using a combination of neural networks and the FEM to accurately solve these types of problems. 
% In our algorithm, the neural network is used to predict the coefficients of the FEM approximation, while the main structure of the numerical scheme is maintained. 
% This allows existing numerical methods such as enriched space methods to be adapted to our network architecture by adding a suitable basis function, as described in references such as \cite{H2020}. 
% This enables our network to effectively capture the behavior of boundary layers, which is an area that has not been well studied in recent machine learning approaches. 

The structure of this paper is outlined as follows: Section \ref{sec:feonet} provides a detailed elaboration on the proposed numerical scheme as well as the theoretical convergence results. Subsequently, Section \ref{sec:exp} presents the experimental results, demonstrating the computational flexibility of FEONet {in both 1D and 2D settings} across various PDE problems. Finally, Section \ref{sec:conclusion} summarizes this work.

\section{Finite element operator network (FEONet)}\label{sec:feonet}
In this section, we present our proposed method, the FEONet, by demonstrating how neural network techniques can be integrated into FEMs to solve the parametric PDEs. Throughout our method description and experiments, we focus on the PDEs of the following form: for uniformly elliptic coefficients $\aaa(\xx)$ and $\varepsilon>0$,
\begin{equation}\label{main_eq}
-\varepsilon\,{\rm{div}}\,(\aaa(\xx)\nabla u)+\mathcal{F}(u)=f\quad{\rm{in}}\,\,D.
\end{equation}
Here $\mathcal{F}$ can be either linear or nonlinear, and both cases will be covered in the experiments performed in Section \ref{sec:exp}. As will be described in more detail later, we shall propose an operator-learning-type method that can provide real-time solution predictions whenever the input data of the PDE changes.

In addition to the explanation for the proposed method, we will also describe the analytic framework of the FEONet, where we can discuss the convergence analysis of the method.

\subsection{Model description}\label{sec:method}
We begin with the following variational formulation of the given PDE \eqref{main_eq}: find $u\in V$ satisfying
\begin{equation}\label{var_for}
B[u,v]:=\varepsilon\int_{D}\aaa(\xx)\nabla u\cdot\nabla v\dx+\int_{D}\mathcal{F}(u)v\dx=\int_{D}fv\dx=:\ell(v)\quad{\text{for all}}\,\,v\in V,
\end{equation}
for a suitable function space $V$. In the theory of FEM, the first step is to define a {\textit{triangulation}} of the given domain $\overline{D}\subset\R^d$, $d=\{1,2,3\}$. We shall assume that $\{\mathcal{T}_h\}_{h>0}$ is a collection of conforming shape-regular triangulations of $\overline{D}\subset\R^d$, consisting of $d$-dimensional simplices $K$. The finite element parameter $h>0$ signifies the maximum mesh size of $\mathcal{T}_h$. In addition, let us denote by $S_h$, the space of all continuous functions $v_h$ defined on $D$ such that the restriction of $v_h$ to an arbitrary triangle is a polynomial. Then we define our finite-dimensional ansatz space as $V_h=S_h\cap V$. We shall denote the set of all vertices in the triangulation by $\{\xx_i\}$, and define the nodal basis $\{\phi_j\}$ for $V_h$, defined by $\phi_j(\xx_i)=\delta_{ij}$.

Then the Galerkin approximation is as follows: find $u_h\in V_h$ satisfying
\begin{equation}\label{gal_approx}
B[u_h,v_h]=\ell(v_h)\quad{\text{for all}}\,\,v_h\in V_h.
\end{equation}
If we write the finite element solution 
\begin{equation}\label{approx_sol}
u_h=\sum^{N(h)}_{k=1}\alpha_k\phi_k,\quad\alpha_i\in\R,
\end{equation}
the Galerkin approximation \eqref{gal_approx} is transformed into the linear algebraic system 
\begin{equation}\label{lin_alg_eq}
A\alpha=F\quad{\text{with}}\,\,A_{ik}:=B[\phi_k,\phi_i]\,\,{\rm{and}}\,\,F_i:=\ell(\phi_i),
\end{equation}
where the matrix $A\in\R^{N(h)\times N(h)}$ is invertible provided that the given PDEs have a suitable structure. Therefore, we can determine the coefficients $\{\alpha_k\}^{N(h)}_{k=1}$ by solving the system of linear equations \eqref{lin_alg_eq}, and consequently, produce the approximate solution $u_h\in V_h$ as defined in \eqref{approx_sol}.

\begin{figure}[t!]
\begin{center}
\includegraphics[width=\textwidth]{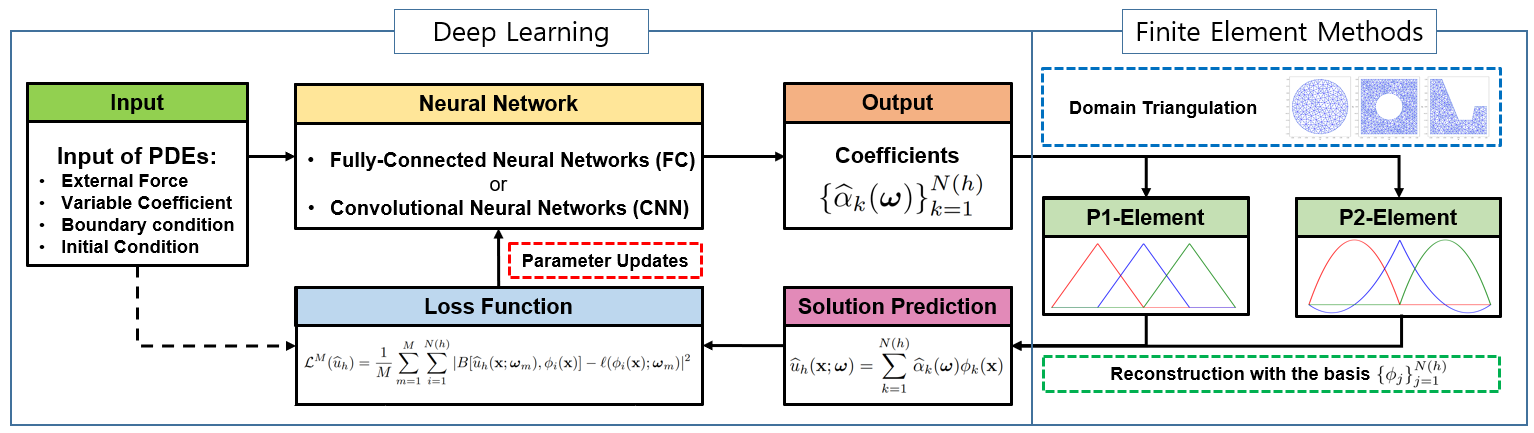}
\end{center}
\caption{Schematic diagram of the Finite Element Operator Network (FEONet).}
\label{fig:scheme}
\end{figure}

Now, in this setting let us describe the method we propose in this paper. In our scheme, an input for the neural network is the data for the given PDEs. In this paper, we choose external force as a prototype input feature. Note, however, that the same scheme can also be developed with other types of data, e.g., boundary conditions, diffusion coefficients, or initial conditions for time-dependent problems. See, for example, Section \ref{subsec:2nd} where extensive experiments were performed addressing various types of input data. For this purpose, the external forcing terms are parametrized by the random parameter $\om$ contained in the (possibly high-dimensional) parameter space $\Omega$. For each realization $f(\xx;\om)$ (and hence for each load vector $F(\om)$ defined in \eqref{lin_alg_eq}), instead of computing the coefficients from the linear algebraic system \eqref{lin_alg_eq}, we approximate the coefficients $\alpha$ by a deep neural network. More precisely, these input features representing external forces pass through the deep neural network and the coefficients $\{\widehat{\alpha}_k\}$ are generated as an output of the neural network. We then reconstruct the solution
\begin{equation}\label{sol_recon}
    \widehat{u}_h(\xx;\om)=\sum_{k=1}^{N(h)}\widehat{\alpha}_k(\om)\phi_k(\xx).
\end{equation}
For the training, we use a residual of the variational formulation \eqref{gal_approx} and define the loss function by summing over all basis functions: for randomly chosen parameters $\om_1,\cdots,\om_M\in\Omega$,
\begin{equation}\label{loss_function}
\LL^M(\widehat{u}_h)=\frac{1}{M}\sum^M_{m=1}\sum^{N(h)}_{i=1}|B[\widehat{u}_h(\xx;\om_m),\phi_i(\xx)]-\ell(\phi_i(\xx);\om_m)|^2 .
\end{equation}

{Note that our proposed loss function is based on the variational form of the given problem, and the parameters are updated in the direction that minimizes the corresponding loss function \eqref{loss_function} defined above. In this paper, we assume that it is possible to find the exact minimizers, allowing us to ignore any errors arising from the optimization process. It is worth noting that loss functions based on the variational form have also been considered in previous work related to PINNs \cite{kharazmi2019variational, KHARAZMI2021113547}. However, unlike our approach, those methods were designed for single-instance predictions and did not utilize the P1 basis elements that are integral to our operator network. Instead, they employed different types of test functions.
While \cite{GOSWAMI2022114587} also proposes an operator-learning method and uses the variational form of governing equations as a loss function along with input-output data pairs, there is a difference from our FEONet approach, which is completely free from the need for input-output data pairs.
The detailed algorithm for training FEONet is provided in the pseudo-code below.}
\floatname{algorithm}{Algorithm}
\begin{algorithm}
{
\caption{{FEONet Training Procedure}}
\label{alg:feonet}
\begin{algorithmic}[1]
\Require Training samples $f(\boldsymbol{x}; \om_m)$ for $m=1, \ldots, M$;
\Ensure The coefficients $\widehat{\alpha}(\om_m)=[\widehat{\alpha}_k(\om_m)]_{k=1, \ldots, N(h)}$ for $m=1, \ldots, M$. Therefore, the predicted solution is $\widehat{u}_h(\boldsymbol{x};\om_m)=\sum_{k=1}^{N(h)}\widehat{\alpha}_k(\om_m)\phi_k(\boldsymbol{x})$;
\State Pre-compute $A$ ($A_{ik}=B[\phi_k,\phi_i]$) and $F(\om_m)$ for $m=1, \ldots, M$ using basis functions $\{\phi_i\}_{i=1}^{N(h)}$ based on domain triangulation $\{\mathcal{T}_h\}$;
\For{$n=1,2,\ldots,N_{\text{epochs}}$}
    \State $\mathcal{L}^M\leftarrow 0$;
    \For{$m=1, \ldots, M$}
            \State Compute $\widehat{\alpha}(\om_m)$ using FEONet with input $f(\boldsymbol{x}; \om_m)$;
            \State $\mathcal{L}^M\leftarrow \mathcal{L}^M + |A\widehat{\alpha}(\om_m)-F(\om_m)|^2$;
    \EndFor
    \State Update FEONet parameters by minimizing $\mathcal{L}^M$;
\EndFor
\end{algorithmic}
}
\end{algorithm}

For each training epoch, once the neural network parameters are updated in the direction of minimizing the loss function \eqref{loss_function}, the external force $f(\xx;\om)$ passes through this updated neural network to generate more refined coefficients, and this procedure is repeated until the sufficiently small loss is achieved. A schematic diagram of the FEONet algorithm is depicted in Figure \ref{fig:scheme}.
\begin{remark}
To avoid confusion for readers, here we provide further clarification: in this paper, the term ``training data'' refers to the labeled data, specifically the pairs of pre-computed PDE data and their corresponding solutions. The term ``input sample'' is used to denote the randomly extracted points from the parametric domain. In this case, since they don't require labeling, there is no computational cost associated with the sampling process. In summary, for the training of FEONet, only input samples are required and the model can be trained without the need for pre-computed training data.
\end{remark}
\begin{remark}
{{There are several reasons why we used the variational residual \eqref{loss_function} as a loss function instead of other types of loss functions, e.g., the physics-informed loss function. Since our method is inherited from the classical FEM, in order to use the well-known techniques from numerical analysis, it is more ``natural'' to consider the loss function we chose. As a consequence, we are able to successfully obtain the theoretical convergence result from the viewpoint of the FEM, where the analysis mainly depends on the form of the loss function. To be more specific, in appendix A.1 and our subsequent paper \cite{feonet_analysis}, we prove the best approximation property of the FEONet, providing the theoretical evidence that we chose the right loss function. In summary, we have adopted the variational residual as a loss function from the theoretical perspective. On the other hand, there are also some advantages of using this weak formulation as a loss function. For example, we can reduce the computation error for the derivative like the variational physics-informed neural networks \cite{kharazmi2019variational, KHARAZMI2021113547}.}}
\end{remark}

\subsection{{Convergence analysis of FEONet}}\label{sect_conv_anal}
{In this section,} we address the convergence result of the FEONet which provides theoretical justification for the proposed method. As outlined in the previous subsection, the proposed method is based on the FEM, and the finite element approximation $u_h$ in \eqref{approx_sol} can be considered as an intermediate solution between the true solution $u$ of \eqref{var_for} and the approximate solution $\widehat{u}_h$ predicted by the FEONet. 

For simplicity, we will focus on self-adjoint equations with homogeneous Dirichlet boundary conditions: 
\begin{equation}\label{self_ad_eq}
    \begin{aligned}
        -{\rm{div}}\, ( \aaa(\xx) \nabla u) + c(\xx) u &= f(\xx)&&{\rm{in}}\,\,D,\\
        u &= 0&&{\rm{on}}\,\,\partial D,
    \end{aligned}
\end{equation}
where the coefficient $\aaa(\xx)$ is uniformly elliptic and $c(\xx)\geq0$.
In the present paper, we focus on the case when the finite element parameter $h>0$ is fixed and deal with the convergence of the FEONet prediction to the finite element solution $u_h$. {The analysis for general equations that cannot be covered by the form \eqref{self_ad_eq} and the theoretical analysis addressing the role of $h>0$ in the error estimate are of independent interest and will be investigated in the forthcoming paper.}

\paragraph{\textbf{Analytic framework}}
As described before, the input of neural networks is an external forcing term $f$, which is parametrized by the random parameter $\om$ in the probability space $(\Omega,\mathcal{G},\mathbb{P})$. In this section, we will regard $f(\xx;\om)$ as a bivariate function defined on $D\times\Omega$, and assume that 
\begin{equation}\label{f_ass}
    f(\xx;\om)\in C(\Omega;L^1(D))\quad{\rm{and}}\quad \sup_{\om\in\Omega}\int_D|f(\xx;\om)|\dx<C
\end{equation}
for some constant $C>0$. For each $\om\in\Omega$, the external force $f(\xx;\om)$ is determined, and the corresponding weak solution is denoted by $u(\xx;\om)$ which satisfies the following weak formulation:
\begin{equation}\label{anal_WF}
    B[u,v]:={\int_D [\aaa(\xx)\nabla u\cdot\nabla v+c(\xx)uv]\dx}=\int_Df(\xx)v\dx=:\ell(v)\quad\forall v\in H^1_0(D).
\end{equation}

For given $h>0$, let $V_h\subset H^1_0(D)$ be a conforming finite element space generated by the basis functions $\{\phi_k\}_{k=1}^{N(h)}$
and $u_h\in V_h$ be a finite element approximation of $u$ satisfying the Galerkin approximation
\begin{equation}\label{anal_Gal}
    B[u_h,v_h]=\ell(v_h)\quad\forall v_h\in V_h.
\end{equation}
We shall write 
\begin{equation}\label{FEM_sol}
    u_h(\xx,\om)=\sum^{N(h)}_{k=1}\alpha^*_k(\om)\phi_k(\xx),
\end{equation}
where $\alpha^*$ is the solution of the linear algebraic equations 
\begin{equation}\label{true_alg}
    A\alpha^*=F,
\end{equation}
with 
\begin{equation}\label{main_mat_vec}
A_{ik}=B[\phi_k,\phi_i]\quad{\rm{and}}\quad F_i=\ell(\phi_i).
\end{equation}
In other words, $u_h$ in \eqref{FEM_sol} is the finite element approximation of the solution to the equation \eqref{self_ad_eq}, as described in Section \ref{sec:method}. It is noteworthy that instead of the equations \eqref{true_alg}, $\alpha^*$ can also be characterized in the following way:
\begin{equation}
    \alpha^*=\argmin_{\alpha\in C(\Omega,\R^{N(h)})}\LL(\alpha),
\end{equation}
where $\LL$ is the population loss
\begin{equation}\label{anal_pop_loss}
    \LL(\alpha)=\mathbb{E}_{\om\sim\mathbb{P}_{\Omega}}\bigg[\sum^{N(h)}_{i=1}|B[\widehat{u}(\om),\phi_i]-\ell(\phi_i;(\om))|^2\bigg]=\|A\alpha(\om)-F(\om)\|^2_{L^2(\Omega)}.
\end{equation}

Next, we shall define a class of neural networks, that plays the role of an approximator for coefficients. For each $L\in\mathbb{N}$, we write an $L$-layer neural network by a function $f^L(x):\R^{n_0}\rightarrow\R^{n_L}$, defined by
\begin{equation}\label{nn_def}
    f^1(x)=W^1x+b^1 \quad{\rm{and}}\quad
    f^{\ell}(x)=W^{\ell}\sigma(f^{\ell-1}(x))+b^{\ell}\,\,\,{\rm{for}}\,\,2\leq\ell\leq L,
\end{equation}
where $W^{\ell}\in\R^{n_{\ell}\times n_{\ell-1}}$ and $b^{\ell}\in\R^{n_{\ell}}$ are the weight and the bias respectively for the $\ell$-th layer, and $\sigma:\R\rightarrow\R$ is a nonlinear activation function. We denote the architecture by the vector ${\vec{\boldsymbol{n}}}=(n_0,\cdots,n_L)$, the class of neural network parameters by $\theta:=\theta_{\vec{\boldsymbol{n}}}=\{(W^1,b^1),\cdots,(W^L,b^L)\}$, and its realization by $\mathcal{R}[\theta](x)$. For a given neural network architecture $\vec{\boldsymbol{n}}$, the family of all possible parameters is denoted by
\begin{equation}\label{nn_para_set}
    \Theta_{\vec{\boldsymbol{n}}}=\left\{\{(W^{\ell},b^{\ell})\}^L_{\ell=1}:W^{\ell}\in\R^{n_{\ell}\times n_{\ell-1}},\,b^{\ell}\in\R^{n_{\ell}}\right\}.
\end{equation}
  For two neural networks $\theta_i$, $i=1,2$, with architectures  $\vec{\boldsymbol{n}}_i=\left(n^{(i)}_0,\cdots,n^{(i)}_{L_i}\right)$, we write $\vec{\boldsymbol{n}}_1\subset\vec{\boldsymbol{n}}_2$ if for any $\theta_1\in\Theta_{\vec{\boldsymbol{n}}_1}$, there exists $\theta_2\in\Theta_{\vec{\boldsymbol{n}}_2}$ satisfying $\mathcal{R}[\theta_1](x)=\mathcal{R}[\theta_2](x)$ for all $x\in\R^{n_0}$.
 
 Now let us assume that there exists a sequence of neural network architectures $\{\vec{\boldsymbol{n}}_n\}_{n\geq1}$ satisfying $\vec{\boldsymbol{n}}_n\subset \vec{\boldsymbol{n}}_{n+1}$ for all $n\in\mathbb{N}$. We define the associated collection of neural networks by
\begin{equation}\label{nn_class}
    \NN^*_n=\{\mathcal{R}[\theta]:\theta\in\Theta_{\vec{\boldsymbol{n}}_n}\}.
\end{equation}
It is straightforward to verify that $\NN^*_n\subset\NN^*_{n+1}$ for any $n\in\mathbb{N}$. We start with the following theorem. This is known as the {\textit{universal approximation theorem}}, and known to hold in various scenarios \cite{UA_1, UA_2, UA_3}.
\begin{theorem}\label{ass_1}
Let $K$ be a compact set in $\R^{d}$ and assume that $g\in C(K,\R^{D})$. Then there holds
\begin{equation}\label{UAT}
    \lim_{n\rightarrow\infty}\inf_{\hat{g}\in\NN^*_n}\|\hat{g}-g\|_{C(K)}=0.
\end{equation}
\end{theorem}

In the present paper, the target function to be approximated by neural networks is $\alpha^*(\om)$. Due to the assumption \eqref{f_ass}, we can verify that $\alpha^*(\om)$ is continuous and bounded. 
This motivates us to deal with the family of bounded neural networks as an ansatz space. For this,
let us write $\sup_{\om\in\Omega}|\alpha^*(\om)|<C_{\alpha}$, and assume that $\sigma^*$ is a bounded activation function with continuous inverse (e.g., sigmoid or tanh activation functions) with $\sigma_{-}:=\inf_{x\in\R}\sigma^*(x)$ and $\sigma_{+}:=\sup_{x\in\R}\sigma^*(x)$. We then define a function $h:[\sigma_{-},\sigma_{+}]\rightarrow[-C_{\alpha},C_{\alpha}]$ satisfying $h(\sigma_{-})=-C_{\alpha}$, $h(\sigma_{+})=C_{\alpha}$. By the assumptions, $h$ also has a continuous inverse. With the above notations,
we consider the following family of neural networks:
\begin{equation}\label{bounded_nn_class}
    \NN_n=\left\{h\circ\sigma^*(g_n):g_n\in\NN_n^*\right\},
\end{equation}
where we can see that
\begin{equation}\label{bddness_NN}
    \|g\|_{C(\Omega)}\leq C_{\alpha},\quad\forall g\in\bigcup_{n\in\mathbb{N}}\NN_n.   
\end{equation}
In the analysis, we will use the following version of the universal approximation theorem for bounded neural networks, which is quoted from \cite{ULGNet_conv}.
\begin{theorem}\label{UAT_main_thm}
Let $\Omega$ be a compact set and $\alpha^*:\Omega\rightarrow\R^{N-1}$ is a continuous function defined in \eqref{true_alg}. Then we have
\begin{equation}\label{UAT_var_3}
    \lim_{n\rightarrow\infty}\inf_{\tilde{g}\in\NN_n}\|\tilde{g}-\alpha^*\|_{C(\Omega)}=0.
\end{equation}
\end{theorem}

Now, with the notation defined above, as a neural network approximation of $\alpha^*$, we seek for $\widehat{\alpha}(n):\Omega\rightarrow\R^{N(h)}$ solving the residual minimization problem
\begin{equation}\label{for_2}
    \widehat{\alpha}(n)=\argmin_{\alpha\in\NN_n}\LL(\alpha),
\end{equation}
where the loss function is minimized over $\NN_n$, and we shall write the corresponding solution by
\begin{equation}\label{for_3}
    u_{h,n}(\xx;\om)=\sum^{N(h)}_{k=1}\widehat{\alpha}(n)_k(\om)\phi_k(\xx).
\end{equation}
Note that for the neural network under consideration $\alpha\in\NN_n$, the input vector is $\om\in\Omega$ which determines the external force $f(\xx;\om)$ and the output is the coefficient vector in $\R^{N(h)}$.

As a final step, let us define the solution to the following discrete residual minimization problem:
\begin{equation}\label{for_4}
    \widehat{\alpha}(n,M)=\argmin_{\alpha\in\NN_n}\LL^M(\alpha),
\end{equation}
where $\LL^M$ is the empirical loss defined by the Monte--Carlo integration of $\LL(\alpha)$:
\begin{equation}\label{anal_emp_loss}
    \LL^M(\alpha)=\frac{|\Omega|}{M}\sum_{m=1}^M\sum^{N(h)}_{i=1}|B[\widehat{u}(\om_m),\phi_i]-\ell(\phi_i;(\om_m))|^2=\frac{|\Omega|}{M}\sum^M_{m=1}|A\alpha(\om_m)-F(\om_m)|^2,
\end{equation}
with $\{\om_n\}^M_{m=1}$ is an i.i.d. sequence of random variables distributed according to $\mathbb{P}_{\Omega}$.
 Then we write the corresponding solution as
\begin{equation}\label{for_5}
    u_{h,n,M}(\xx;\om)=\sum^{N(h)}_{k=1}\widehat{\alpha}(n,M)_k(\om)\phi_k(\xx),
\end{equation}
which is the numerical solution actually computed by the scheme proposed in the present paper. 

In order to provide appropriate theoretical backgrounds for the proposed method, it would be reasonable to show that our solution is sufficiently close to the approximate solution computed by the proposed scheme for various external forces, as the index $n\in\mathbb{N}$ for a neural-network architecture and the number of input samples $M\in\mathbb{N}$ goes to infinity. Mathematically, it can be formulated as
\begin{equation}\label{main_error}
    \|u-u_{h,n,M}\|_{L^2(\Omega;L^2(D))}\rightarrow0\quad{\rm{as}}\,\,h\rightarrow0\,\,{\rm{and}}\,\,n,M\rightarrow\infty.
\end{equation}
The main error can be split into three parts:
\begin{equation}\label{error_split}
    u-u_{h,n,M}=
    (u-u_{h})+
    (u_{h}-u_{h,n})+
    (u_{h,n}-u_{h,n,M}).
\end{equation}
The first error is the one caused by the finite element approximation, which is assumed to be negligible for a proper choice of $h$. 
In fact, according to the classical theory of the FEM (e.g. \cite{BS08, C02}) if $h>0$ is sufficiently small, then we can expect to obtain a sufficiently small error. Therefore, in this paper, we shall assume that we have chosen suitable $h>0$ which guarantees a sufficiently small error for the finite element approximation, and only focus on the convergence of $u_h-u_{h,n,M}$ with a fixed $h>0$. 

The second error is referred to as the approximation error, as it appears when we approximate the
target function with a class of neural networks. The third error is often called the generalization error, which measures how well our approximation predicts solutions for unseen data. We will focus on proving that our approximate solution converges to the finite element solution (which is sufficiently close to the true solution) as a neural network architecture grows and the number of input samples goes to infinity, i.e., we will show that
\[
u_h-u_{h,n,M}\rightarrow0\quad{\rm{as}}\,\,n,M\rightarrow\infty.
\]

\paragraph{\textbf{Approximation error}}
We first note from \eqref{anal_pop_loss} and \eqref{anal_emp_loss}, that $A$ defined in \eqref{FEM_sol}, \eqref{true_alg} mostly determines the structures of the loss functions and it would be advantageous for us to analyze the loss function if we know more about the matrix. Note that the matrix $A$ is determined by the structure of the given differential equations, the choice of basis functions, and the boundary conditions. Therefore, the characterization of $A$ which is useful for the analysis of the loss function and can cover a wide range of PDE settings simultaneously is important. The next lemma addresses this issue, which is quoted from \cite{ULGNet_conv}.

\begin{lemma}\label{matrix_thm}
    Suppose that the matrix $A\in\R^{N(h)\times N(h)}$ is symmetric and invertible, and let $\rho_{\min}=\min_i\{|\lambda_i|\}$, $\rho_{\max}=\max_i\{|\lambda_i|\}$ where $\{\lambda_i\}$ is the family of eigenvalues of the matrix $A$. Then we have for all $\xx\in\R^{N(h)}$,
    \begin{equation}\label{eigen_est}
        \rho_{\min}|\xx|\leq|A\xx|\leq\rho_{\max}|\xx|.
    \end{equation}
\end{lemma}
Since the equation under consideration \eqref{self_ad_eq} is self-adjoint, the corresponding bilinear form $B[\cdot,\cdot]$ defined in \eqref{anal_WF} is symmetric which automatically guarantees that the matrix $A$ in our case is symmetric. Furthermore, due to the facts that the coefficient $\aaa(\cdot)$ is uniformly elliptic and $c(\cdot)$ is non-negative, the bilinear form $B[\cdot,\cdot]$ is coercive, which implies that $A$ is positive-definite. Therefore, we can apply Lemma \ref{matrix_thm} to the matrix $A$ of our interest.

We are now ready to present our first convergence result regarding the approximation error. The proof is based on Lemma \ref{matrix_thm} and can be found in Appendix \ref{app_1}.
\begin{theorem}\label{approx_main_thm}
Assume that \eqref{f_ass} holds. Then we have that
\begin{equation}\label{app_t_est}
\|\alpha^*-\widehat{\alpha}(n)\|_{L^2(\Omega)}\rightarrow0\quad{\rm{as}}\,\,n\rightarrow\infty.
\end{equation}
\end{theorem}

\paragraph{\textbf{Generalization error}}
In order to handle the generalization error, we first introduce the concept so-called {\textit{Rademacher complexity}} \cite{Rade_3}. 
\begin{definition}
For a collection $\{X_i\}_{i=1}^M$ of i.i.d. random variables, we define the Rademacher complexity of the function class $\mathcal{G}$ by
\[
    R_M(\mathcal{G})=\mathbb{E}_{\{X_i,\varepsilon_i\}^M_{i=1}}\bigg[\sup_{f\in\mathcal{G}}\bigg|\frac{1}{M}\sum^M_{i=1}\varepsilon_if(X_i)\bigg|\bigg],
\]
where $\varepsilon_i$'s denote i.i.d. Bernoulli random variables, in other words, $\mathbb{P}(\varepsilon_i=1)=\mathbb{P}(\varepsilon_i=-1)=\frac{1}{2}$ for all $i=1,\cdots,M$.
\end{definition}

Note that the Rademacher complexity is the expectation of the maximum correlation between
the random noise $(\varepsilon_1,\cdots,\varepsilon_M)$ and the vector $(f(X_1),\cdots,f(X_M))$, where the supremum is taken over the class of functions $\mathcal{G}$. By intuition, the Rademacher complexity of $\mathcal{G}$ is often used to measure how the functions from $\mathcal{G}$ can fit random noise. For a more detailed discussion on the Rademacher complexity, see \cite{Rade_1, Rade_2, Rade_4}.

Now we define the function class of interest
\begin{equation}\label{fct_class}
    \mathcal{G}_n:=\{|A\alpha-F|^2:\alpha\in\NN_n\},
\end{equation}
where the matrix $A$ and the vector $F$ are defined in \eqref{main_mat_vec}. In the following theorem, note that we assume that for any $n\in\mathbb{N}$, the Rademacher complexity of $\mathcal{G}_n$ converges to zero, which is known to be true in many cases \cite{Rade_4, rade_upper_2, rade_upper_3}. For the proof of the theorem below, see Appendix \ref{app_2}.%\cite{rade_upper_1, rade_upper_2,rade_upper_3,rade_upper_4}.

\begin{theorem}\label{gen_conv_thm}
Suppose that \eqref{f_ass} holds, and we further assume that for any $n\in\mathbb{N}$, $\lim_{M\rightarrow\infty}R_M(\mathcal{G}_n)=0$. Then we have with probability $1$ that
\[
    \lim_{n\rightarrow\infty}\lim_{M\rightarrow\infty}\|\widehat{\alpha}(n,M)-\widehat{\alpha}(n)\|_{L^2(\Omega)}=0.
\]
\end{theorem}

\paragraph{\textbf{Convergence of FEONet}}
From Theorem \ref{approx_main_thm} and Theorem \ref{gen_conv_thm} combined with the triangular inequality, we have probability 1 over i.i.d.
samples that
\begin{equation}\label{conv_1}
    \lim_{n\rightarrow\infty}
    \lim_{M\rightarrow\infty}\|\alpha^*-\widehat{\alpha}(n,M)\|_{L^2(\Omega)}=0.
\end{equation}
Based on the convergence of the approximate coefficients \eqref{conv_1}, let us now state the convergence theorem for the predicted solution, whose proof is presented in Appendix \ref{app_3}.
\begin{theorem}\label{main_thm_whole}
Suppose \eqref{f_ass} holds, and assume that for all $n\in\mathbb{N}$, $R_M(\widetilde{{\mathcal{G}}}_n)\rightarrow0$ as $M\rightarrow \infty$, with $\widetilde{{\mathcal{G}}}_n:=\{|A\alpha-F|^2:\alpha\in\NN_n\}$. Then for given $h>0$, we have with probability $1$ over i.i.d. samples that
\begin{equation}\label{main_conv_whole}
    \lim_{n\rightarrow \infty}\lim_{M\rightarrow \infty}\|u_h-u_{h,n,M}\|_{L^2(\Omega;H^s(D))}=0.
\end{equation}
\end{theorem}

\section{Numerical experiments}\label{sec:exp}
In this section, we demonstrate the experimental results of the FEONet across various settings for PDE problems. We randomly generate 3000 input samples for external forces and train the FEONet in an unsupervised manner without the usage of pre-computed input-output $(f,u)$ pairs. We then evaluate the performance of our model with another randomly generated test set consisting of 3000 pairs of $f$ and the corresponding solutions $u$. For this test data, we employed the FEM with a sufficiently fine mesh discretization to ensure that the numerical solutions can be considered as true solutions (See Appendix \ref{appendix:subsec:input}). While there are other operator learning-based methods, such as PINO \cite{PINO}, that can be trained without input-output data pairs, they often face limitations when dealing with domains of complex shapes.
% A detailed comparison with the results for PINO can be found in Appendix \ref{appendix:subsec:compare_PINO_PIDON}.
{Our primary comparison will be against DeepONet \cite{lu2021learning} and POD-DeepONet \cite{lu2022comprehensive}, with varying numbers of training data, as well as against PIDeepONet \cite{PIDON}. The FEONet is designed to learn the operator without requiring pre-computed training data, which significantly reduces computational overhead. The generation of random input samples, which are used for unsupervised learning, imposes negligible computational costs, offering flexibility in choosing the number of input function samples for training. In contrast, DeepONet (and POD-DeepONet) relies on supervised learning, requiring input-output data pairs. This necessitates pre-computing the output data for each input function sample, resulting in substantial computational demands. Consequently, it is challenging to directly compare the training times of FEONet with those of DeepONet (or POD-DeepONet). Note that training FEONet to solve the 2D Poisson equation in a domain with a circular hole, achieving a 1.3\% relative error, took approximately 28 minutes and 36 seconds on an NVIDIA RTX A5000 24GB GPU.} The triangulation of domains and the generation of the nodal basis for the FEM are performed using FEniCS \cite{AlnaesEtal2015,LoggEtal2012}. {A detailed description of the model used in the experiments, as well as the comparison of computational time between FEM and operator learning models, can be found in Appendix \ref{appendix:subsec:hyperparam} and \ref{appendix:subsec:compare_time}.}

In Section 3.1, we will solve PDE problems in various environments and settings to verify the computational flexibility of our method and confirm the flexibility of the proposed operator learning method by employing various types of input functions. In Section 3.2, we will conduct some further experiments on the parabolic problem and the Stokes problem. In Section 3.3, we will demonstrate that the FEONet framework can make an accurate solution prediction for the singular perturbation problems. Finally, in Section 3.4, additional experiments will be conducted to confirm the robustness and reliability of our method across various performance metrics. Throughout the experiments, we will cover both one-dimensional and two-dimensional cases, and we have appropriately organized them according to the objectives of each experiment.

\begin{table}[t!]
  \caption{Mean Rel. $L^2$ test errors $\mathbf{(\times10^{-2})}$ with standard deviations $\mathbf{(\times10^{-2})}$ for the 3000 test set on diverse PDE problems. Five training trials are performed independently. Domain \RNum{1} (circle), Domain \RNum{2} (square with a hole), Domain \RNum{3} (polygon) with Eq. \eqref{eq_domain}, BC \RNum{1} (Dirichlet), BC \RNum{2} (Neumann) with Eq. \eqref{eq_bdry}, Eq \RNum{1} (second-order linear equation \eqref{eq_varcoeff}) and Eq \RNum{2} (Burgers' equation \eqref{eq_burgers}).}
  \label{table:total_error}
  \centering
  \resizebox{\columnwidth}{!}{
  \begin{tabular}{cccccccc}
    \toprule
    Model (\#Train data) & Domain \RNum{1} & Domain \RNum{2} & Domain \RNum{3} & BC \RNum{1} & BC \RNum{2} & Eq \RNum{1} & Eq \RNum{2}\\
    \midrule
    DON (supervised, w/30) & 27.15\scriptsize{$\pm$1.16} & 51.21\scriptsize{$\pm$3.58} & 53.92\scriptsize{$\pm$4.59} & 21.75\scriptsize{$\pm$1.19} & 22.75\scriptsize{$\pm$1.05} & 24.38\scriptsize{$\pm$1.37} & 10.26\scriptsize{$\pm$0.14} \\
    DON (supervised, w/300) & 2.10\scriptsize{$\pm$0.75} & 5.62\scriptsize{$\pm$0.37} & 6.22\scriptsize{$\pm$0.96} &
    0.68\scriptsize{$\pm$0.11} & 0.96\scriptsize{$\pm$0.06} & 0.76\scriptsize{$\pm$0.10} & \textbf{0.20}\scriptsize{$\pm$0.09} \\
    DON (supervised, w/3000) & \textbf{0.69}\scriptsize{$\pm$0.17}  & 4.75\scriptsize{$\pm$0.75} & 6.20\scriptsize{$\pm$1.00} & 0.53\scriptsize{$\pm$0.36} & 0.33\scriptsize{$\pm$0.09} & 0.33\scriptsize{$\pm$0.27} & 0.24\scriptsize{$\pm$0.13} \\
    PIDeepONet (w/o labeled data) & 9.80\scriptsize{$\pm9.41$} & 101.03\scriptsize{$\pm$167.46} & 32.89\scriptsize{$\pm6.34$} & 1.51\scriptsize{$\pm0.46$} & 1.43\scriptsize{$\pm0.45$} & 19.41\scriptsize{$\pm11.30$} & 2.66\scriptsize{$\pm0.71$} \\
    Ours (w/o labeled data) & 1.24\scriptsize{$\pm$0.00} & \textbf{1.76}\scriptsize{$\pm$0.03} & \textbf{0.51}\scriptsize{$\pm$0.00} & \textbf{0.13}\scriptsize{$\pm$0.01} & \textbf{0.32}\scriptsize{$\pm$0.03} & \textbf{0.13}\scriptsize{$\pm$0.01} & 0.54\scriptsize{$\pm$0.07} \\
    \bottomrule
  \end{tabular}
  }
\end{table}

\begin{table}[t!]
{
  \caption{Mean Rel. $L^2$ test errors $\mathbf{(\times10^{-2})}$ for the 1000 test set on diverse PDE problems. Coefficient (Variable coefficient as an input) with Eq. \eqref{eq_varcoeff}, Boundary (Boundary condition as an input) with Eq. \eqref{eq_bdry_input}, Initial (Time-dependent problem) with Eq. \eqref{eq_time} and Vector-valued (Stokes problem) with Eq. \eqref{eq_stokes}.}
  \label{table:total_error2}
  \centering
  \resizebox{\columnwidth}{!}{
  \begin{tabular}{ccccc}
    \toprule
    Model (\#Train data) & Coefficient & Boundary & Initial  & Vector-valued\\
    \midrule
    POD-DeepONet (supervised, w/10) & 25.276 & 47.500 & 7446.077 & - \\
    POD-DeepONet (supervised, w/100) & 10.536 & 28.471 & 230.288 & - \\
    POD-DeepONet (supervised, w/1000) & \textbf{0.371} & 6.095 & 10.628 & - \\
    Ours (w/o labeled data) & 0.612 & \textbf{0.430} & \textbf{0.218} & \textbf{5.141} \\
    \bottomrule
  \end{tabular}
  }
  }
\end{table}

\subsection{Performance evaluation of FEONet in various scenarios}\label{subsec:1st}
To evaluate the computational flexibility of the FEONet, we conducted a sequence of experiments covering diverse domains, boundary conditions, and equations.

\begin{figure}[htbp]
\begin{center}
\includegraphics[width=0.85\textwidth]{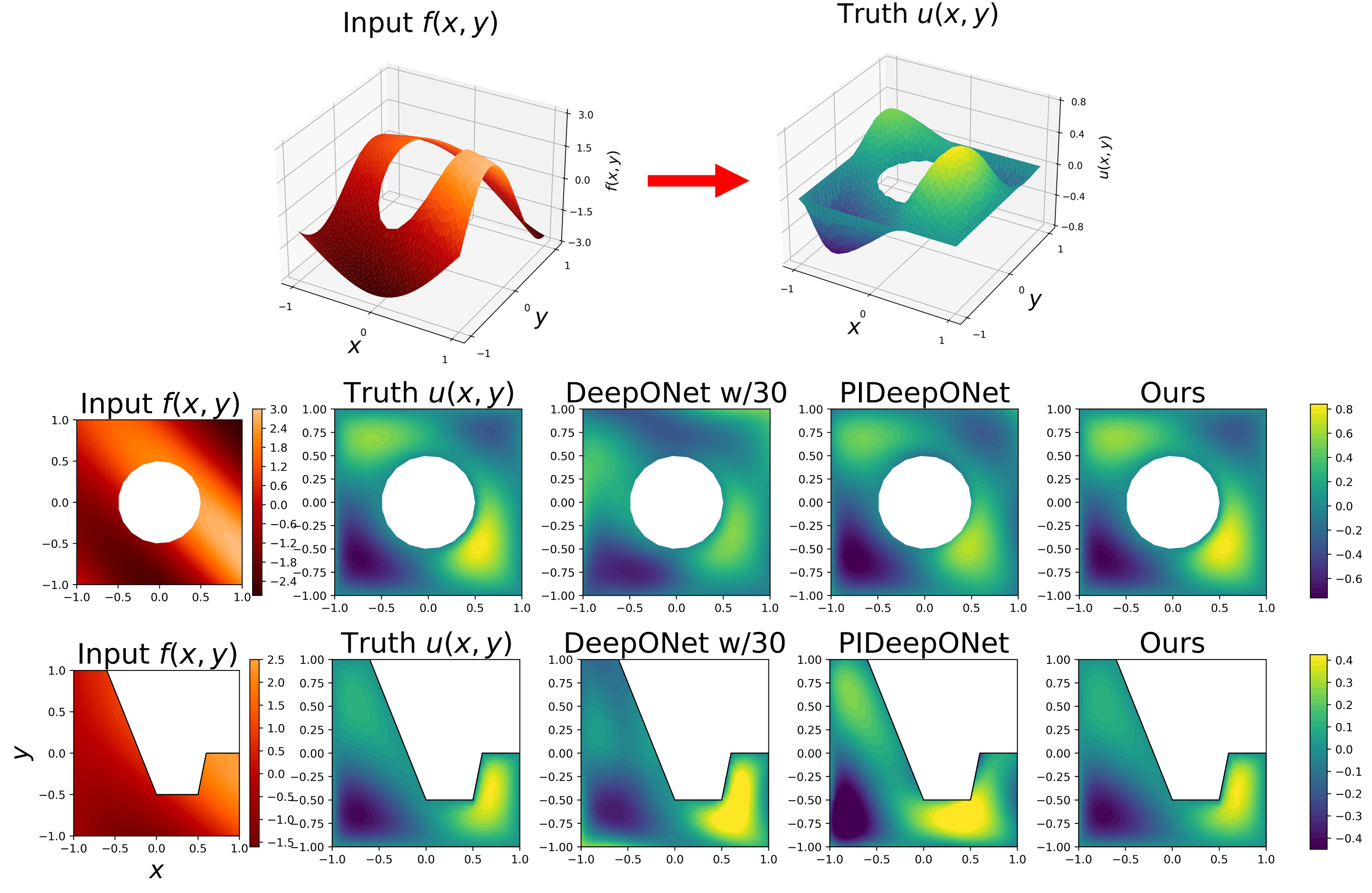}
\end{center}
\caption{Solution profiles for complex geometries.}
\label{fig:domain_profile}
\end{figure}

\subsubsection{Simulations on different domains}
We consider the 2D convection-diffusion equation as
\begin{equation}\label{eq_domain}
    \begin{aligned}
        -\varepsilon \Delta u+\boldsymbol{v}\cdot \nabla u&=f(x,y), && (x,y) \in D,\\
        u(x,y)&=0, && (x,y) \in \partial D.
    \end{aligned}
\end{equation}
For our numerical experiments, we fixed the values of $\varepsilon=0.1$ and $\boldsymbol{v}=(-1,0)$ across various domains $D$, which included a circle, a square with a hole, and a polygon (see Figure \ref{fig:domain}). 
The second to fourth columns of Table \ref{table:total_error} display the mean relative $L^2$ errors of the test set for FEONet (w/o labeled data), PIDeepONet (w/o labeled data), and DeepONet (with 30/300/3000 input-output data pairs) across various domains. When the DeepONet is trained with either 30 or 300 data pairs, FEONet consistently surpasses them in terms of numerical errors. Notably, even when DeepONet is trained with 3000 data pairs, FEONet achieves errors that are either comparable to or lower than those of DeepONet.
It's noteworthy that in more intricate domains like Domain \RNum{2} and Domain \RNum{3}, both the DeepONet (even with 3000 data pairs) and PIDeepONet find it challenging to produce accurate predictions. Figure \ref{fig:domain_profile} demonstrates the solution profile, emphasizing that FEONet is more precise in predicting the solution than both DeepONet and PIDeepONet, particularly in complex domains.

\subsubsection{Simulations on Dirichlet and Neumann boundary conditions}
We tested the performance of each model when the given PDEs were equipped with either Dirichlet or Neumann boundary conditions. We consider the convection-diffusion equation as
\begin{equation}\label{eq_bdry}
  \begin{aligned}
    -\varepsilon u_{xx}+bu_x&=f(x), && x \in D,\\
    u(x)=0\quad\text{or}\quad &u_x(x)=0, && x \in \partial D,
  \end{aligned} 
\end{equation}
where $\varepsilon=0.1$ and $b=-1$. For the Neumann boundary condition, we add $u$ to the left-hand side of \eqref{eq_bdry} to address the uniqueness of the solution. The FEONet has an additional significant advantage to make the predicted solutions satisfy the exact boundary condition, 
by selecting the appropriate coefficients for the basis functions. Using the FEONet without any input-output training data, we are able to obtain similar or slightly higher accuracy compared to the DeepONet with 3000 training data, for both homogeneous Dirichlet or Neumann boundary conditions as shown in fifth and sixth columns in Table \ref{table:total_error}.

\subsubsection{Simulations on linear and nonlinear PDEs}
We performed some additional experiments applying the FEONet to the various equations: (i) general second-order linear PDE with variable coefficients defined as
\begin{equation}\label{eq_varcoeff}
  \begin{aligned}
    -\varepsilon u_{xx}+b(x)u_x+c(x)u=f(x),\;\;\;\;\;&x\in D,\\
    u(x)=0,\;\;\;\;\;&x\in \partial D,
  \end{aligned} 
\end{equation}
where $\varepsilon=0.1$, $b(x)=x^2+1$ and $c(x)=x$; and (ii) the Burgers' equation defined as
\begin{equation}\label{eq_burgers}
  \begin{aligned}
    -u_{xx}+uu_x&=f(x), && x \in D,\\
    u(x)&=0, && x \in \partial D.
  \end{aligned} 
\end{equation}
Since the Burgers' equation is nonlinear, an iterative method is required for classical numerical schemes including the FEM, which usually causes some computational costs. Furthermore, it is more difficult than linear equations to obtain training data for nonlinear equations. One important advantage of the FEONet is its ability to effectively learn the nonlinear structure without any training data. Once the training process is completed, the model can provide accurate real-time solution predictions without relying on iterative schemes.
% See Appendix \ref{appendix:subsec:burgers} for a detailed explanation of the loss function for the nonlinear Burgers' equation. 
The seventh and eighth columns of Table \ref{table:total_error} highlight the effectiveness of the FEONet, predicting the solutions with reasonably low errors for various equations. 
Although the eighth column of Table \ref{table:total_error} shows that the DeepONet with 3000 training data has a lower error compared to our model, it demonstrates the strength of the FEONet that achieves a similar level of accuracy even for the nonlinear PDE, that can be trained in an unsupervised manner.

\begin{figure}[t!]
\begin{center}
\includegraphics[width=\textwidth]{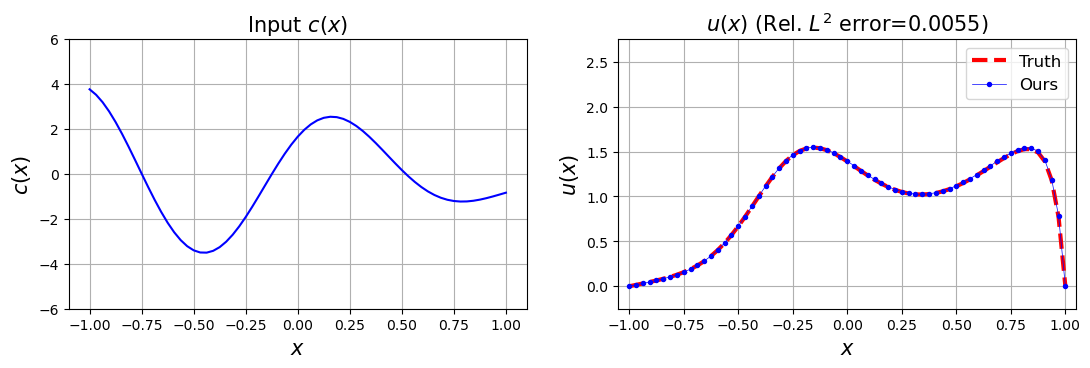}
\end{center}
\caption{Solution profiles predicted by the trained FEONet for PDE \eqref{eq_varcoeff} obtained from variable coefficients $c(x)$ as an input function.}
\label{fig:input_varcoeff}
\end{figure}

\subsubsection{Simulations on various types of input functions}\label{subsec:2nd}
In this section, we will demonstrate that our method can accommodate various types of input functions. Specifically, we shall utilize, as an input for our model, variable coefficients and boundary conditions. Throughout the series of experiments, we will confirm that our approach can make a fast and accurate solution prediction for such various scenarios. Note, however, that the proposed method cannot take the geometry of the domain as an input of neural networks, which is of independent interest (See, for example, the references \cite{domain_geo_1, domain_geo_2} where this issue is addressed).

\paragraph{\textbf{Variable coefficient as an input}}
The FEONet is capable of learning the operator mapping $\mathcal{G}:c(x)\mapsto u(x)$ from the variable coefficient $c(x)$ to the corresponding solution $u(x)$ of a PDE \eqref{eq_varcoeff}. Figure \ref{fig:input_varcoeff} displays the solution profile obtained when variable coefficients are considered as input. {The second column in Table \ref{table:total_error2} shows that when the input is a coefficient variable, POD-DeepONet, using all 1,000 data pairs, produces a slightly smaller error than FEONet, but the errors are quite similar. This demonstrates the strength of FEONet in achieving comparable accuracy while being trained in an unsupervised manner.}
% This highlights the flexibility of FEONet in adapting to a wide range of input function types, thereby underscoring its versatility.

\begin{figure}[]
\begin{center}
\includegraphics[width=\textwidth]{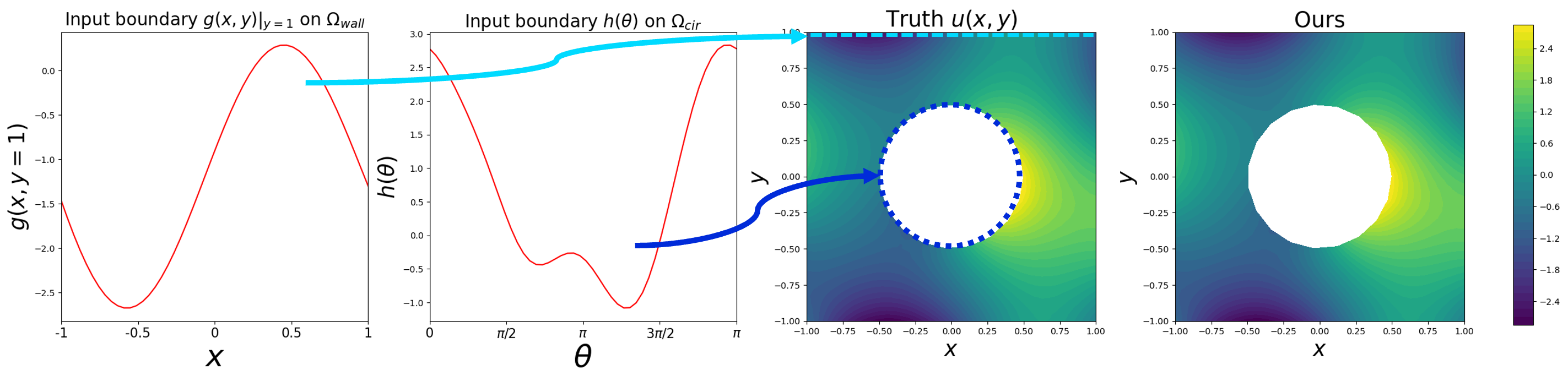}
\end{center}
\caption{Solution profiles predicted by the trained FEONet for PDE \eqref{eq_bdry_input} obtained from the boundary condition input. Note that the boundary condition $g(x,y)$ for the region with a circle hole inside $\Omega_{\text{cir}}$ is represented as $h(\theta)$ in polar coordinates. %{The relative $L^2$ error, when comparing the exact solution profile with the predicted solution for this specific test case, is 0.0074. Meanwhile, the mean relative $L^2$ test error across the entire test set is 0.0043.}\
The average of relative $L^2$ errors across the entire test set when comparing the exact solution to the predictions is 0.0043.
}
\label{fig:input_bdry}
\end{figure}

\paragraph{\textbf{{Boundary condition as an input}}}
We address the case when the input is a boundary condition. We consider the 2D Poisson equation as
\begin{equation}\label{eq_bdry_input}
    \begin{aligned}
        - \Delta u(x,y)&=f(x,y), && (x,y) \in D,\\
        u(x,y)&=g(x,y), && (x,y) \in \partial D.
    \end{aligned}
\end{equation}
where $f(x,y)=2$ and the domain $D$ is a square with a hole (see (b) in Figure \ref{fig:domain}). The FEONet is utilized to learn the operator $\mathcal{G}:g(x,y)\mapsto u(x,y)$. Figure \ref{fig:input_bdry} demonstrates that even in the case of an operator that takes the boundary condition $g(x,y)$ as input, the FEONet can predict the solution with high accuracy, {as shown in the third column of Table \ref{table:total_error2}, compared to POD-DeepONet.} This demonstrates the versatility and adaptability of the FEONet approach for accommodating different types of input functions.

\subsection{Further extensions of FEONet}
In this section, we will extend FEONet to handle the parabolic problem and the Stokes problem.
\subsubsection{Time-dependent problem}

\begin{figure}[h!]
\begin{center}
\includegraphics[width=\textwidth]{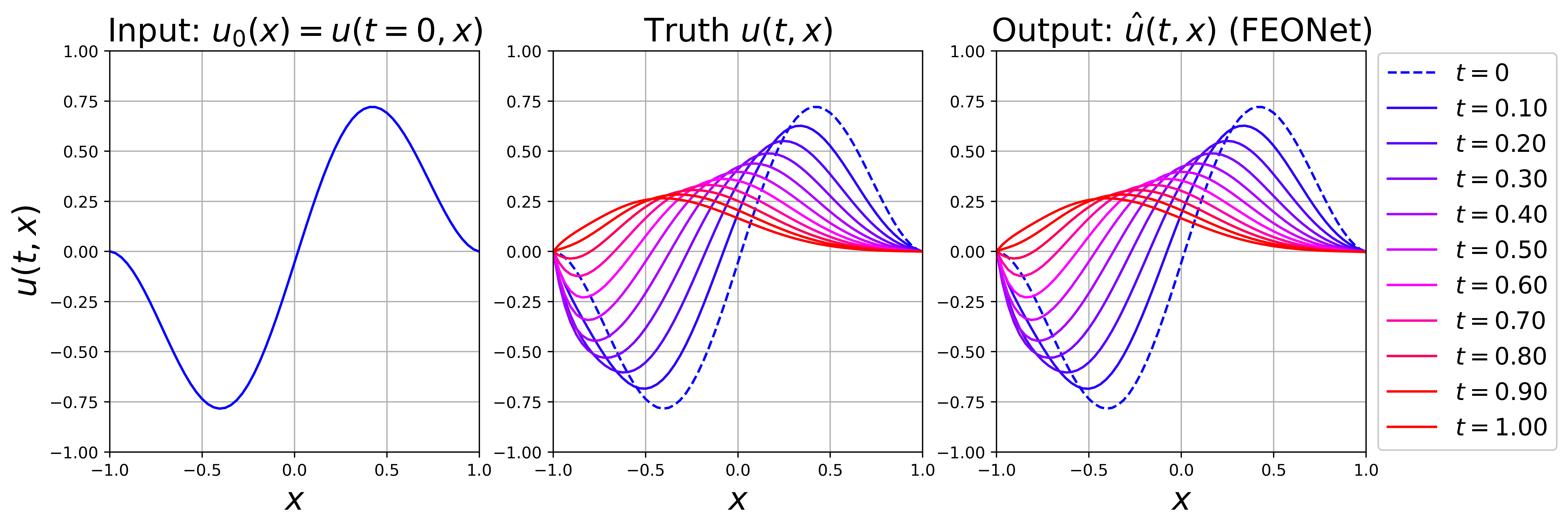}
\end{center}
\caption{Solution profiles predicted by the trained FEONet for PDE \eqref{eq_time} obtained from the initial condition $u(t=0,x)$ as an input function. 
%{The relative $L^2$ error, when comparing the exact solution profile with the predicted solution for this specific test case, is 0.0028. Meanwhile, the mean relative $L^2$ test error across the entire test set is 0.0010.}
The average of relative $L^2$ error across the entire test set, when comparing the exact solution profile with the predicted solution, is 0.0022.
}
\label{fig:input_initial}
\end{figure}

{The FEONet can be extended to handle the parabolic problem with initial conditions as input.} We consider the time-dependent problem with varying initial conditions. We consider the 1D convection-diffusion equation
\begin{equation}\label{eq_time}
  \begin{aligned}
    u_t-\nu u_{xx}+bu_x&=0, && t\in [0,1], x \in D=[-1,1],\\
    u(0,x)&=u_0(x), && x \in D,\\
    u(t,x)&=0, && x \in \partial D,
  \end{aligned} 
\end{equation}
where $\nu=0.1$ and $b=-1$. We aim to learn the operator $\mathcal{G}:u_0(x)\mapsto u(t,x)$. We use the implicit Euler method to discretize the time with $\Delta t=0.01$. We employed the marching-in-time scheme proposed by \cite{krishnapriyan2021characterizing} to sequentially learn and predict solutions over the time domain $t=[0,1]$, divided into 10 intervals. For training FEONet, we utilized only the input function $u_0(x)$ and equation \eqref{eq_time} without any additional data. {See Appendix \ref{appendix:subsec:time_d_pde} for more details.} Figure \ref{fig:input_initial} illustrates FEONet's effective prediction of true solutions using an initial condition from the test dataset. {The fourth column in Table \ref{table:total_error2} shows that the FEONet—trained without any data pairs—yields a smaller error than POD-DeepONet.}

\subsubsection{Stokes problem}

\begin{figure}[]
\begin{center}
\includegraphics[width=\textwidth]{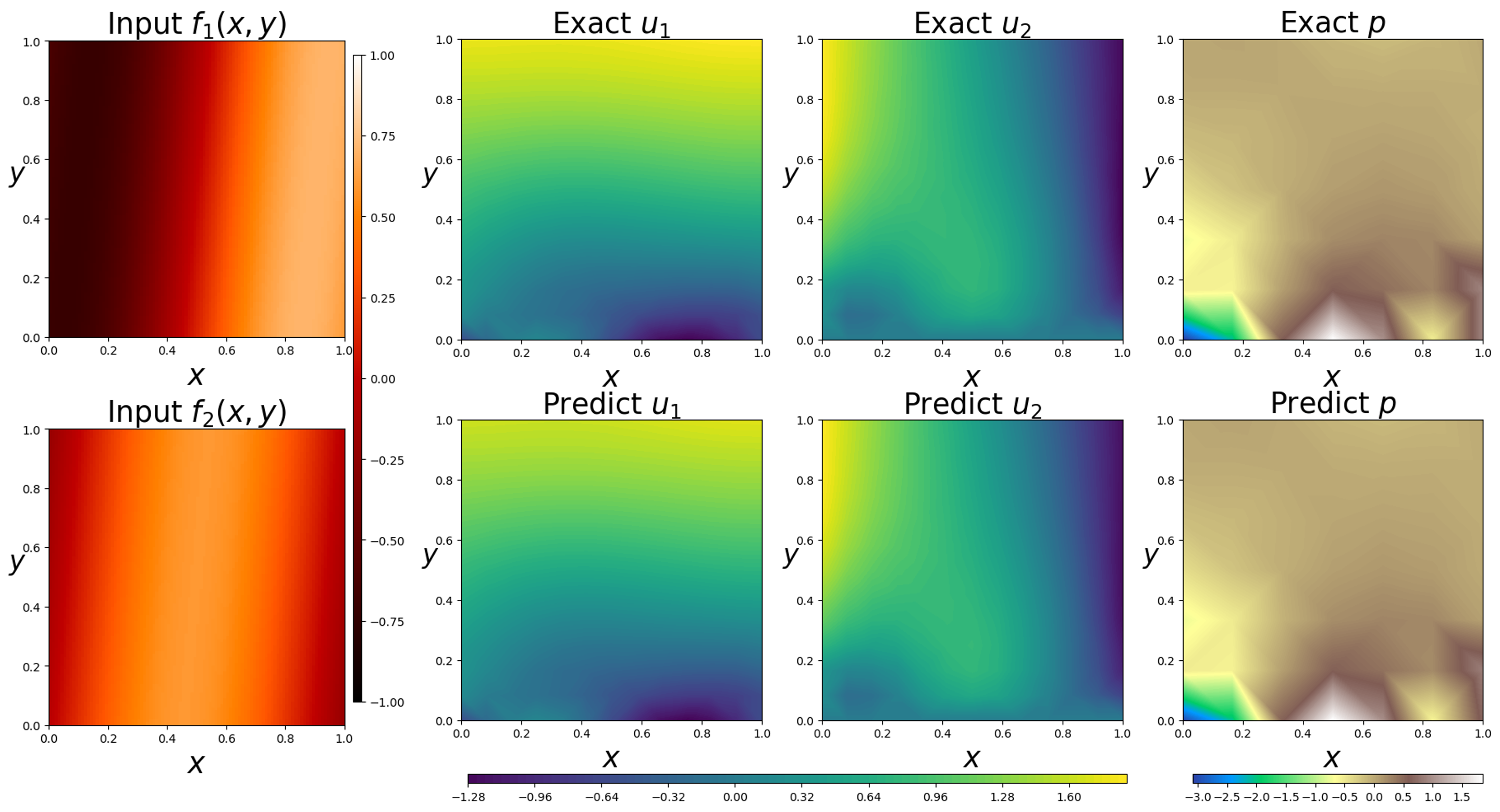}
\end{center}
\caption{Solution profiles predicted by the trained FEONet for 2D Stokes equation \eqref{eq_stokes}. 
%For this specific test case, the relative $L^2$ errors when comparing the exact solution to the predictions for $u_1$, $u_2$, and $p$ are 0.0421, 0.0395, and 0.0580, respectively. Meanwhile, the average relative $L^2$ errors for $u_1$, $u_2$, and $p$ across the entire test set are 0.0499, 0.1092, and 0.0685.
The average of relative $L^2$ errors across the entire test set when comparing the exact solution to the predictions for $u_1$, $u_2$, and $p$ are 0.0499, 0.1092, and 0.0685, respectively.
}
\label{fig:stokes}
\end{figure}
We consider a more intricate yet practical fluid dynamics issue: the two-dimensional steady-state Stokes equations, which can be described as follows:
\begin{equation}\label{eq_stokes}
    \begin{aligned}
        - \nabla\cdot(\nabla \boldsymbol{u}-pI)&=\boldsymbol{f}, && (x,y) \in D=[0,1]^2,\\
        \nabla\cdot \boldsymbol{u}&=0, && (x,y) \in D,\\
        \boldsymbol{u}&=\boldsymbol{g}, && (x,y) \in \Gamma_D=\{(x,y)\in D|y=0\},\\
        (\nabla \boldsymbol{u}-pI)\cdot \boldsymbol{n}&=(0,0), && (x,y) \in \Gamma_N=\{(x,y)\in D|x=0,1\text{ or }y=1\}
    \end{aligned}
\end{equation}
where $\boldsymbol{g}(x,y)=(3+1.7\sin(2\pi x),0)$. 
%In this context, the velocity vector $\boldsymbol{u}=(u_1,u_2)$ and the pressure $p$ are adopted from  \cite{11122/8316}. 
We aim to learn the operator $\mathcal{G}:\boldsymbol{f}=(f_1,f_2)\mapsto[\boldsymbol{u},p]$ using the FEONet. 
% {\color{magenta}{The variational form of \eqref{eq_stokes} is written as:
% \begin{equation}
% \begin{split}
% a((u, p), (v, q))
%                         &= \int_{D} \nabla u \cdot \nabla v
%          - \nabla \cdot v \ p
%          + \nabla \cdot u \ q \, {\rm d} x, \\
% L((v, q))
%                         &= \int_{D} f \cdot v \, {\rm d} x.
% \end{split}
% \end{equation}
% }}
For this, we employ the P2 Lagrange element for velocity and the P1 Lagrange element for pressure, which is known as the Taylor-Hood element \cite{elman2014finite}.
Figure \ref{fig:stokes} demonstrates the results of the operator learning for the 2D Stokes equation using FEONet. This illustrates that FEONet is applicable to the learning of operators for the system of PDEs with vectorial inputs and outputs. {In this problem, which requires predicting three outputs—$u_1(x,y)$, $u_2(x,y)$, and $p(x,y)$—the experiment for POD-DeepONet was omitted as shown in the fifth column of Table \ref{table:total_error2} because constructing a natural extension of the basis using POD is challenging. This remains an area for future work involving the application of POD-DeepONet.}

\subsection{Singular perturbation problem}\label{subsec:singular}
\begin{figure}[htbp]
\begin{center}
\includegraphics[width=\textwidth]{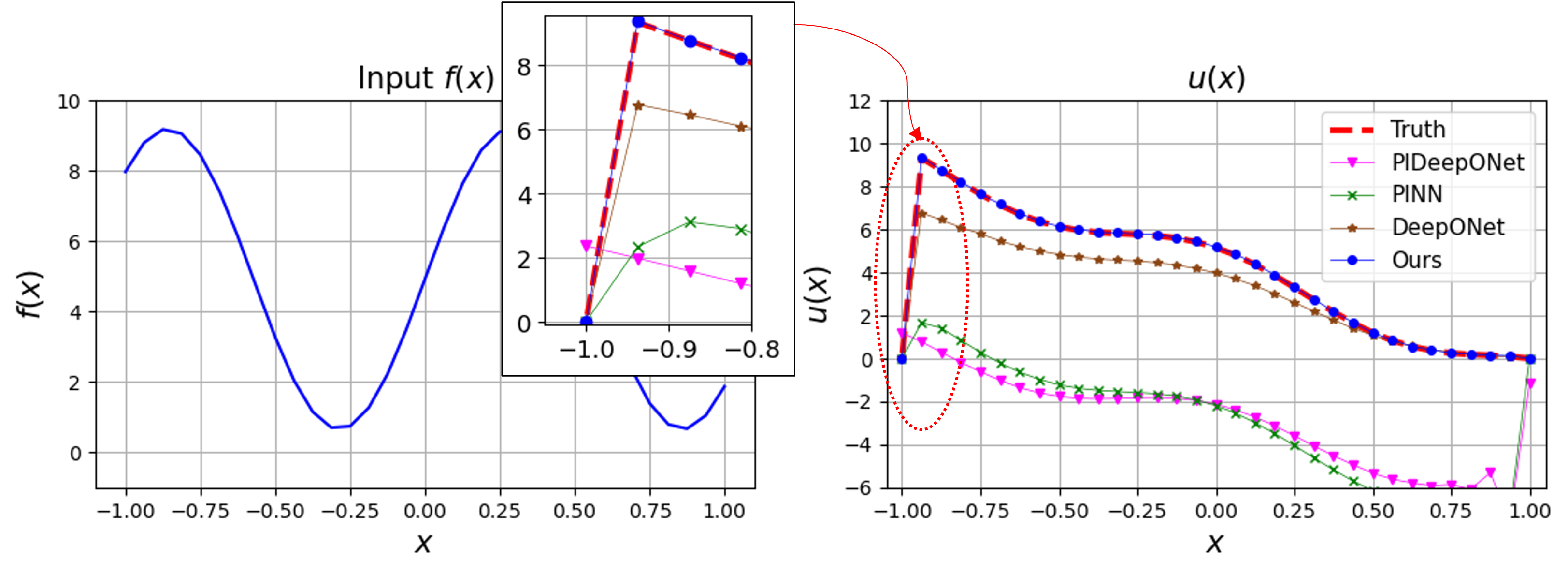}
\end{center}
\caption{Input function $f$ and the corresponding approximate solution $u$ obtained by various deep-learning-based methods for the singular perturbation problem with $\varepsilon=10^{-5}$.}
\label{fig:bdrylayer}
\end{figure}

{Since the proposed method is based on domain meshing, we may face some potential challenges associated with the scenarios involving singular behavior including interior or boundary shocks. In particular, as our method is not mesh-free, refined meshes are required to handle certain problems, such as shock formation. This is a recognized limitation of our approach. However, while our method may face difficulties in solving singular problems, it benefits from the ability to leverage various existing numerical analysis techniques to address these challenges. To the best of our knowledge, such problems are very challenging for other machine learning-based approaches, including mesh-free methods.} More precisely, one of the notable advantages of the FEONet, which predicts coefficients through well-defined basis functions, is its applicability to problems involving singularly perturbed PDEs. 
The boundary layer problem is a well-established challenge in scientific computing. This issue arises from the presence of thin regions near boundaries where the solution exhibits steep gradients. Such regions significantly impact the overall behavior of the system.
Typically, in cases where the target function demonstrates abrupt transitions, such as when the diffusion coefficient is small in our model problem, neural network algorithms frequently encounter difficulties in converging to optimal solutions. This is largely due to the phenomenon known as spectral bias \cite{rahaman2019spectral}. Neural networks are generally predicated on a smooth prior assumption, but spectral bias leads to a failure to accurately capture abrupt transitions or singular behaviors of the target solution function.
To address this challenge, we propose a novel semi-analytic machine learning approach, guided by singular perturbation analysis, for effectively capturing the behavior of thin boundary layers.
This is accomplished by incorporating additional basis functions guided by theoretical considerations.
For this we consider
\begin{equation}\label{eq_singular}
  \begin{aligned}
    -\varepsilon u_{xx}+bu_x&=f(x), && x \in D,\\
    u(x)&=0, && x \in \partial D,
  \end{aligned} 
\end{equation}
where $\varepsilon \ll 1$. For the implementation, we assign the values of $b=-1$ and $\varepsilon=10^{-5}$.
For instance, in the right panel of Figure \ref{fig:bdrylayer}, the red dotted line (ground-truth) shows the solution profile which produces a sharp transition near $x=-1$. 
To capture the boundary layer, we utilize an additional basis function known as the corrector function in mathematical analysis, defined as
\begin{equation}\label{eq:corrector}
    \phi_{\text{cor}}(x):=e^{-(1+x)/\varepsilon}-(1-(1-e^{-2/\varepsilon})(x+1)/2).
\end{equation}

\begin{table}[htbp]
\caption{Comparisons of various deep-learning-based models for solving the singular perturbation problem. Each model is trained five times independently.}
  \label{table:bdrylayer_error}
  \centering
  \resizebox{\columnwidth}{!}{
\begin{tabular}{llllll}
\hline
Model             & Ours & PINN & DeepONet w/30 & PIDeepONet \\ \hline
Requirement for labeled data      & \textbf{No}  & \textbf{No}  & Yes  & \textbf{No}    \\
Multiple instance & \textbf{Yes}  & No  & \textbf{Yes} & \textbf{Yes}    \\
Mean Rel. $L^2$ test error\scriptsize{$\pm$std}   &  \textbf{0.0132}\scriptsize{$\pm$0.0091}    &  1.3827\scriptsize{$\pm$0.7580} &  0.2303\scriptsize{$\pm$0.0074}   &     0.5713\scriptsize{$\pm$0.0007}   \\ 
 Meshless & No & Yes & Yes & Yes\\
\hline
\end{tabular}
} 
\end{table}

% see Appendix \ref{appendix:subsec:singular} for more detailed information regarding the derivation of the additional basis function. 
By integrating the boundary layer element into the finite element space, we establish the proposed enriched Galerkin space for utilization in the FEONet. As the corrector basis is incorporated alongside the conventional nodal basis functions in the FEM, we also predict the additional coefficient originating from the corrector basis. Note that the computational cost remains minimal as the enriched basis exclusively encompasses boundary elements.
Table \ref{table:bdrylayer_error} summarizes the characteristics of various operator learning models and the corresponding errors when solving the singularly perturbed PDE described in \eqref{eq_singular}.
The performance of the FEONet in accurately solving boundary layer problems surpasses that of other models, even without labeled data, across multiple instances. %Due to the inclusion of a very small $\varepsilon$ in the PDE's residual loss, operator learning with the PINN and the PIDeepONet is difficult or unstable.
The incorporation of a very small $\varepsilon$ in the residual loss of the PDE makes operator learning challenging or unstable when using PINN and PIDeepONet.
%In fact, neural networks have a smooth prior, which causes them to suffer from the boundary layer problem and results in unstable training due to the direct calculation of PDE residual. 
In fact, neural networks inherently possess a smooth prior, which can lead to difficulties in handling boundary layer problems. 
In contrast, FEONet utilizes theory-guided basis functions, enabling the predicted solution to accurately capture the sharp transition near the boundary layer.
The PIDeepONet faces difficulties in effectively learning both the residual loss and boundary conditions, while PINNs struggle to predict the solution.
Figure \ref{fig:bdrylayer} exhibits the results of operator learning for the singular perturbation problem.
As depicted in the zoomed-in graph of Figure \ref{fig:bdrylayer}, the FEONet stands out as the only model capable of accurately capturing sharp transitions, whereas other models exhibit noticeable errors.

\subsection{Generalization capacity and convergence rate}\label{subsec:further}
\paragraph{\textbf{Generalization capacity}}
The FEONet is a model capable of learning the operator without the need for paired input-output training data. What we need for the training is the generation of random input samples, which only imposes negligible computational cost. Therefore, the number of input function samples for the FEONet training can be chosen arbitrarily. In the experiments performed in Section \ref{sec:exp}, we generated 3000 random input function samples to train the FEONet. On the other hand, the DeepONet requires paired input-output training data for supervised learning which causes significant computational costs to prepare. In this paper, we have focused our experiments on considering the DeepONet trained with 30, 300, and 3000 input-output data pairs. This allowed us to highlight the FEONet's ability to achieve a certain level of accuracy without relying on data pairs.

\begin{figure}[htbp]
\begin{center}
\includegraphics[width=\textwidth]{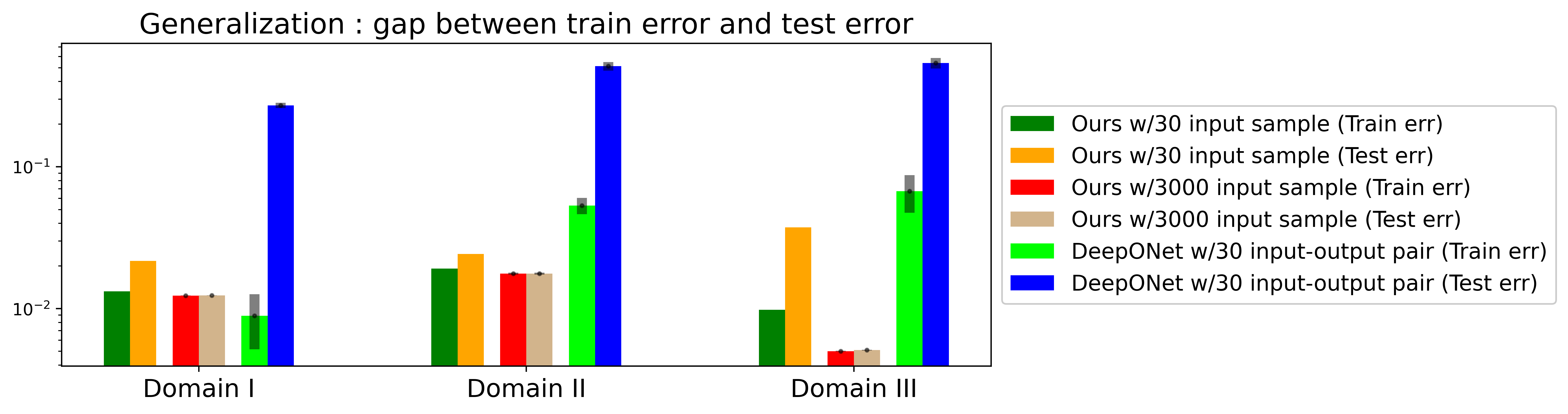}
\end{center}
\caption{Rel. $L^2$ train and test error comparison between the FEONet with 30/3000 input samples (w/o labeled data) and the DeepONet (supervised) with 30 train-test data pair.}
\label{fig:generalization2}
\end{figure}

{Figure \ref{fig:generalization2} presents the results comparing the generalization capability of FEONet when trained with only 30 random input samples (green and yellow bars), which should be distinguished from the 30 input-output data pairs used for supervised learning.
More specifically, this figure illustrates the generalization capabilities of FEONet in comparison to DeepONet across the three distinct domain settings discussed in Section 3.1: Domain I (circle), Domain II (square with a hole), and Domain III (polygon). Each bar represents the models' performance when solving Equation (3.1) within these different geometrical domains.
Figure \ref{fig:generalization2} compares the training and test errors of FEONet, trained through an unsupervised approach with 30 and 3,000 input samples, against DeepONet, trained through supervised learning with 30 input-output data pairs. The difference between the training and test errors serves as a reliable indicator of the models' generalization effectiveness. Notably, FEONet demonstrates superior generalization compared to DeepONet, even without requiring input-output paired data. Remarkably, with only 30 input samples for training, FEONet achieves significantly better generalization and exhibits lower test errors than DeepONet. The experimental details can be found in Appendix \ref{appendix:subsec:generalization}.}

\begin{figure}[htbp]
\begin{center}
\includegraphics[width=\textwidth]{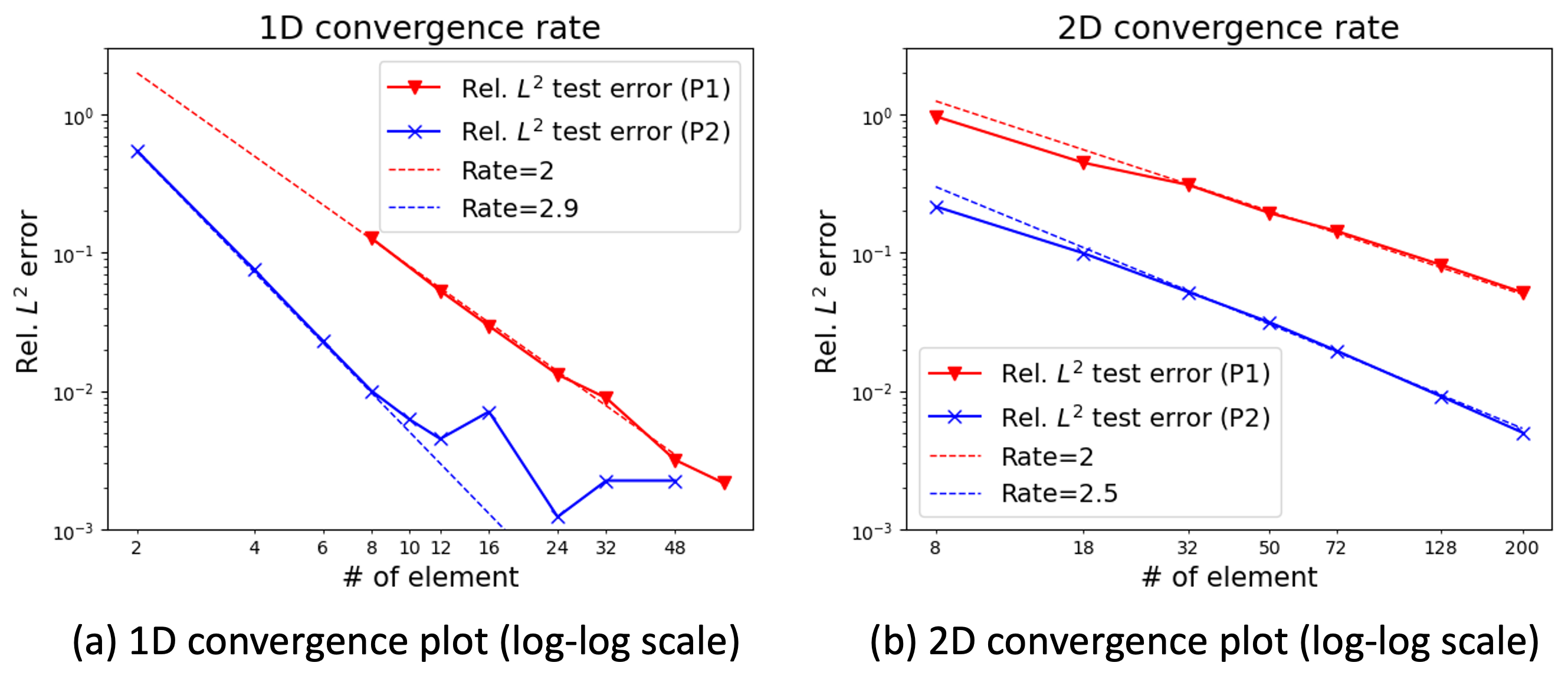}
\end{center}
\caption{Rel. $L^2$ error of the FEONet with respect to the number of elements.}
\label{fig:rate}
\end{figure}

\paragraph{\textbf{Rate of convergence}}
% \begin{wrapfigure}{r}{0.35\textwidth}
% \vspace{-10mm}
% \hspace{0mm}
% \includegraphics[width=0.35\textwidth]{rate.png}
% \caption{Rel. $L^2$ error of the FEONet with respect to the number of elements.}
% \label{fig:rate}
% \vspace{-4mm}
% \end{wrapfigure}
One of the intriguing aspects is that the theoretical results on the convergence error rate of the FEONet are also observable in the experimental results. Figure \ref{fig:rate} shows the relationship between the test error and the number of elements, utilizing both P1 and P2 basis functions. In 1D problems, the convergence rate is approximately 2 for P1 basis functions and 2.9 for P2 basis functions. Meanwhile, in 2D problems, the convergence rate is approximately 2 for P1 basis functions and 2.5 for P2 basis functions.  These remarkable trends confirm the theoretical convergence rates observed in the experimental results. The experimental details can be found in Appendix \ref{appendix:subsec:convergence}.

{
\begin{remark}
    As we can see from the left in Figure \ref{fig:rate}, as we increased the number of elements, i.e., on a high-fidelity mesh, the error began to increase again after a certain point, which means that there is an issue regarding the high-dimensional coefficient
output of the neural network. To explain this phenomenon, we theoretically investigated the FEONet in our subsequent paper \cite {feonet_analysis} and successfully identified the cause of such phenomenon. More precisely, we discovered that the condition number of the finite element matrix is closely related to the high-fidelity issue of FEONet, which makes total error increase when the degree of freedom is large. Furthermore, this theoretical finding provides a critical clue for resolving the computational issue arising from high-fidelity meshes. By using a preconditioner to define a loss function, we can significantly reduce the total error and the FEONet can effectively make a solution prediction even for the case of a large number of elements. See \cite {feonet_analysis} for details.
\end{remark}
}

\section{Conclusions}\label{sec:conclusion}
%PDEs are prevalent in various problems, including fluid mechanics, electromagnetics, and geophysics, and play a prominent role in the mathematical analysis, modeling, and simulation of complex physical systems. 
In this paper, we introduce the FEONet, a novel deep learning framework that approximates nonlinear operators in infinite-dimensional Banach spaces using finite element approximations.
Our operator learning method for solving PDEs is both simple and remarkably effective, resulting in significant improvements in predictive accuracy, domain flexibility, and data efficiency compared to existing techniques.
Furthermore, we demonstrate that the FEONet is capable of learning the solution operator of parametric PDEs, even in the absence of paired input-output training data, and accurately predicting solutions that exhibit singular behavior in thin boundary layers.
The FEONet is capable of learning the operator by incorporating not only external forces but also boundary conditions, variable coefficients, and the initial conditions of time-dependent PDEs. {On the other hand, our method also has some drawbacks that need to be overcome. For example, since the proposed method is not mesh-free and based on domain meshing, we may face challenges in solving high-dimensional PDEs or handling singular problems, such as shock formation (see, e.g. \cite{PINN007}). Note, however, that since our method is based on the classical FEM, we can directly apply numerical analysis techniques to solve these problems, which will be addressed in future work.} If these aspects are effectively supplemented, we are confident that the FEONet can provide a new and promising approach to simulate and model intricate, nonlinear, and multiscale physical systems, with a wide range of potential applications in science and engineering.

\appendix
\section{Proof of theoretical results}\label{appen_conv}

\subsection{Approximation error}\label{app_1}
We shall prove Theorem $2.4$ addressing the approximation error for the coefficients.
\begin{proof}[Proof of Theorem \ref{approx_main_thm}]
Since $A$ is symmetric and invertible, by Proposition \ref{matrix_thm}, we obtain
\begin{align*}
    \|\alpha^*-\widehat{\alpha}(n)\|^2_2
    &\lesssim\|A\alpha^*-A\widehat{\alpha}(n)\|^2_2\lesssim\left(\|A\alpha^*-F\|^2_2+\|A\widehat{\alpha}(n)-F\|^2_2\right)\\
    &=\LL(\widehat{\alpha}(n))\lesssim\inf_{\alpha\in\NN_n}\LL(\alpha)=\inf_{\alpha\in\NN_n}\|A\alpha-F\|^2_2\\
    &\lesssim\inf_{\alpha\in\NN_n}\left(\|A\alpha-A\alpha^*\|^2_2+\|A\alpha^*-F\|^2_2\right)\\
    &\lesssim\inf_{\alpha\in\NN_n}\|\alpha-\alpha^*\|^2_2.
\end{align*}
Note that the constants appearing in the inequalities above may depend on $h>0$. But as we mentioned before, we assumed that $h>0$ is fixed with the sufficiently small finite element error.
Then by the universal approximation property , $\inf_{\alpha\in\NN_n}\|\alpha-\alpha^*\|^2_2\rightarrow0$ as $n\rightarrow\infty$, which completes the proof. 
\end{proof}

\subsection{Generalization error}\label{app_2}
Next, we shall find the connection between the generalization error and the Rademacher complexity for the uniformly bounded function class $\mathcal{G}$. In the following theorem, we assume that the function class is $b$-uniformly bounded, meaning that $\|f\|_{\infty}\leq b$ for arbitrary $f\in\mathcal{G}$.

\begin{theorem}\label{rade_thm}
{\rm{[Theorem 4.10 in \cite{Rade_2}]}} Assume that the function class $\mathcal{G}$ is $b$-uniformly bounded and let $M\in\mathbb{N}$. Then for arbitrary small $\delta>0$, we have
\[
\sup_{f\in\mathcal{G}}\bigg|\frac{1}{M}\sum^M_{i=1}f(X_i)-\mathbb{E}[f(X)]\bigg|\leq 2R_M(\mathcal{G})+\delta,
\]
with probability at least $1-\exp(-\frac{M\delta^2}{2b^2})$.
\end{theorem}

Next, note from Lemma \ref{matrix_thm} that
\[
    \|A\alpha-F\|_{L^{\infty}(\Omega)}\leq\|A\alpha\|_{L^{\infty}(\Omega)}+\|F\|_{L^{\infty}(\Omega)}\lesssim\rho_{\max}\|\alpha\|_{L^{\infty}(\Omega)}+\|f\|_{C(\Omega;L^1(D))}.
\]
Since the class of neural networks under consideration is uniformly bounded and \eqref{f_ass} holds, we see that for any $n\in\mathbb{N}$, $\mathcal{G}_n$ is $\tilde{b}$-uniformly bounded for some $\tilde{b}>0$.
Then the following lemma is the direct consequence of Theorem \ref{rade_thm} within our setting.
\begin{lemma}\label{rade_main}
Let $\{\om_m\}^M_{m=1}$ be i.i.d. random samples selected from the distribution $\mathbb{P}_{\Omega}$. Then for arbitrary small $\delta>0$, we obtain with probability at least $1-2 \exp(-\frac{M\delta^2}{32 \tilde{b}^2})$ that
\begin{equation}\label{loss_diff}
    \sup_{\alpha\in\NN_n}\left|\LL^M(\alpha)-\LL(\alpha)\right|\leq2R_n(\mathcal{G}_n)+\frac{\delta}{2}.
\end{equation}
\end{lemma}

Now from Lemma \ref{rade_main}, we shall prove Theorem $2.6$ concerning the generalization error for the coefficients. 
\begin{proof}[Proof of Theorem \ref{gen_conv_thm}]
By Proposition \ref{matrix_thm} and the definition \eqref{for_2}, we obtain
\begin{equation}
    \begin{aligned}\label{gen_pr_mid}
        \|\widehat{\alpha}(n)-\widehat{\alpha}(n,M)\|^2_2
        &\lesssim\|A\widehat{\alpha}(n)-A\widehat{\alpha}(n,M)\|^2_2\lesssim\left(\|A\widehat{\alpha}(n)-F\|^2_2 + \|A\widehat{\alpha}(n,M)-F\|^2_2\right)\\
        &=\left(\LL(\widehat{\alpha}(n))+\LL(\widehat{\alpha}(n,M))\right)\lesssim\LL(\widehat{\alpha}(n,M)).
    \end{aligned}
\end{equation}
We next apply Lemma \ref{rade_main} for $\delta=2M^{-\frac{1}{2}+\varepsilon}$ with $0<\varepsilon<\frac{1}{2}$. Then with probability at least $1-2\exp(-\frac{M^{2\varepsilon}}{8\tilde{b}^2})$, we have from the minimality of $\widehat{\alpha}(n,M)$ that,
\[
    \LL(\widehat{\alpha}(n,M))
    \leq \LL^M(\widehat{\alpha}(n,M))+2R_M(\mathcal{G}_n)+M^{-\frac{1}{2}+\varepsilon}\leq \LL^M(\widehat{\alpha}(n))+2R_M(\mathcal{G}_n)+M^{-\frac{1}{2}+\varepsilon}.
\]
Using Lemma \ref{rade_main} again gives us that
\[
    \LL(\widehat{\alpha}(n,M))\leq \LL(\widehat{\alpha}(n))+4R_M(\mathcal{G}_n)+2M^{-\frac{1}{2}+\varepsilon}.
\]
Letting $M\rightarrow\infty$ on \eqref{gen_pr_mid}, we obtain that
\[
    \lim_{M\rightarrow\infty}\|\widehat{\alpha}(n,M)-\widehat{\alpha}(n)\|^2_2\lesssim\LL(\widehat{\alpha}(n)).
\]
As we did before, we conclude that
\begin{align*}
    \lim_{n\rightarrow\infty}\lim_{M\rightarrow\infty}\|\widehat{\alpha}(n,M)-\widehat{\alpha}(n)\|^2_2
    &\lesssim \lim_{n\rightarrow\infty}\LL(\widehat{\alpha}(n))\lesssim \lim_{n\rightarrow\infty}\inf_{\alpha\in\NN_n}\|\alpha-\alpha^*\|^2_2= 0, 
\end{align*}
which completes the proof.
\end{proof}

\subsection{Convergence of FEONet}\label{app_3}
Based on the results proved above, now we can address the convergence of the predicted solution by the FEONet to the finite element solution.
\begin{proof} [Proof of Theorem \ref{main_thm_whole}]
From Theorem \ref{approx_main_thm}, Theorem \ref{gen_conv_thm}, we note that for a fixed $h>0$,
\begin{align*}
    &\|u_h-u_{h,n,M}\|^2_{L^2(\Omega;L^2(D))}\\
    &=\int_{\Omega}\int_D\bigg|\sum^{N(h)}_{i=1}(\alpha^*_i-\widehat{\alpha}(n,M)_i)\phi_i\bigg|^2\dx\dom\\
    &=\int_{\Omega}\int_D\bigg|\sum^{N(h)}_{i,j=1}(\alpha^*_i-\widehat{\alpha}(n,M)_i)(\alpha^*_j-\widehat{\alpha}(n,M)_j)\phi_i\phi_j\bigg|^2\dx\dom\\
    &\leq\int_{\Omega}\int_D\bigg(\sum^{N(h)}_{i,j=1}|\alpha^*_i-\widehat{\alpha}(n,M)_i|^2|\phi_i|^2+\sum^{N(h)}_{i,j=1}|\alpha^*_j-\widehat{\alpha}(n,M)_j|^2|\phi_j|^2\bigg)\dx\dom\\
    &\leq\int_{\Omega}\int_D2N(h)\sum^{N(h)}_{k=1}|\alpha^*_k-\widehat{\alpha}(n,M)_k|^2|\phi_k|^2\dx\dom\\
    &\leq 2|D|N(h)\|\alpha^*-\widehat{\alpha}(n,M)\|^2_{L^2(\Omega)}\rightarrow0,
\end{align*}
where we have used the fact that $|\phi_j|\leq1$ for all $1\leq j\leq N(h)$ together with Young's inequality. 
\end{proof}

\section{Experimental details}\label{appendix:sec:experiment}
\subsection{Types and random generation of input functions for FEONet}\label{appendix:subsec:input}
% This paper specifically has focused on the learning of the solution operator of PDEs from external forcing terms.
%  We shall confirm this based on some additional experiments with various input functions.
In order to train the network, we generate random input functions (e.g. external forcing functions, variable coefficient functions, boundary conditions, and initial conditions). Inspired by \cite{bar2019learning}, we created a random signal $f(\xx;\om)$ as a linear combination of sine functions and cosine functions. More precisely, we use
\begin{equation}
    f(x)=m_0\sin(n_0x)+m_1\cos(n_1x)
\end{equation}for a 1D case and 
\begin{equation}
    f(x,y)=m_0\sin(n_0x+n_1y)+m_1\cos(n_2x+n_3y)
\end{equation} for a 2D case where $m_i$ for $i=1,2$ and $n_j$ for $j=0,1,2,3$ are drawn independently from the uniform distributions.
% whose ranges are shown in the second column of Table \ref{table:parameters}.
In our experiments for DeepONet with input-output data pairs, we prepare a total of 30/300/3000 pairs of $(f,u)$ by applying the FEM on sufficiently fine meshes.
It is worth noting that even when considering different random input functions, such as those generated by Gaussian random fields, we consistently observe similar results. This robustness indicates the reliability and stability of the FEONet approach across various input scenarios.

{
\subsection{Details on hyperparameters}\label{appendix:subsec:hyperparam}
For the problems under consideration in the paper, we used neural networks, which consist of 6 convolutional layers with swish activation followed by a fully connected layer flattening the output. For 1D problems, we used Conv1D, while Conv2D was used for 2D problems. The FEONet was trained with the LBFGS optimizer along with the following hyperparameters.
\begin{itemize}
    \item Maximal number of iterations per optimization step: $10$,
    \item Termination tolerance on first-order optimality : $10^{-15}$,
    \item Termination tolerance on function value/parameter  changes: $10^{-15}$,
    \item Update history size: $10$.
\end{itemize}
For both the DeepONet (or POD-DeepONet) and the PIDeepONet, a fully connected neural network with a depth of 3 and a width of 100 was used for the trunk and branch networks. The two models differ only in the loss function, where the training data and the residual of the PDE were used. The PINN also used a fully connected neural network with a depth of 3 and a width of 100. All these models used the tanh activation function and a learning rate of $10^{-4}$ with the commonly used ADAM optimizer. To ensure the sufficient convergence of the training results for the boundary layer problem, we conducted five experiments with the PINN for $5\times10^5$ iterations and five experiments with the PIDeepONet for $10^4$ iterations. Note that computing the residual loss of PDE for every input function $f$ in the PIDeepONet requires significant computation time using the automatic differentiation.
}

{
\subsection{Comparison of computational time}\label{appendix:subsec:compare_time}
Once the operator-learning ML models (DeepONet, POD-DeepONet, and FEONet) are trained, they can quickly infer solutions for varying PDE parameters. This capability is a key reason why operator networks are becoming increasingly popular in scientific machine learning. In contrast, solving PDEs with traditional FEM for various parameters requires re-solving the PDE each time the parameters change, incurring computational costs. This challenge is further exacerbated when dealing with nonlinear problems, where iterative methods must be applied for each set of PDE parameters, adding to the computation overhead. As mentioned earlier, FEONet is particularly well-suited to address this issue, highlighting a critical difference between traditional FEM and FEONet. To quantify this advantage, let us consider a scenario where 10,000 parametric PDEs with varying parameters need to be solved—a common requirement in computational material database generation. In such cases, operator learning methods like FEONet are significantly faster than traditional FEM. Table \ref{table:compare_fem} illustrates the inference time comparison based on the number of input samples, underscoring the efficiency of models like FEONet in fast-inference applications.}
\begin{table}[h!]
{
\centering
\caption{Comparison of computational times for FEM, DeepONet, POD-DeepONet and FEONet with different numbers of input samples.}
\label{table:compare_fem}
\centering
\begin{tabular}{cccccc}
\hline
\# input samples & 1 & 10 & 100 & 1000 & 10000 \\ \hline
FEM & 0.012 (s) & 0.141 (s) & 1.039 (s) & 9.979 (s) & 99.826 (s) \\
DeepONet & 0.147 (s) & 0.123 (s) & 0.251 (s) & 0.338 (s) & 0.746 (s) \\ 
POD-DeepONet & 0.170 (s) & 0.190 (s) & 0.191 (s) & 0.398 (s) & 0.765 (s) \\
FEONet (Ours) & 0.106 (s) & 0.221 (s) & 0.315 (s) & 0.374 (s) & 1.529 (s) \\ \hline
\end{tabular}
}
\end{table}

{
\subsection{Extension of FEONet to time-dependent PDEs}\label{appendix:subsec:time_d_pde}
We use the FEONet to learn an operator that maps $u_0(x)=u(t=0,x)$ to $u(\Delta t,x)$ (more precisely, the coefficients of $u(\Delta t,x)$) in Equation \eqref{eq_time}. The other model then predicts the solution at the next time step solution $u(2\Delta t,x)$ as an input $u(\Delta t,x)$),  and continues iteratively. To facilitate this, we employed the implicit Euler method to discretize time with $\Delta t=0.01$, and specifically trained the model using the following loss function:
\begin{multline}\label{eq:loss_time-d}
    \mathcal{L}^M(\widehat{u}_h)=\frac{1}{M}\sum^M_{m=1}\sum^{N}_{i=1}\bigg| \int_D \frac{\widehat{u}_h(l\Delta t,x)-\widehat{u}_h((l-1)\Delta t,x)}{\Delta t}\phi_i\mathrm{d}x\\
+ \nu \int_D \widehat{u}_h(l\Delta t,x)_x(\phi_i)_x\mathrm{d}x
- b\int_D\widehat{u}_h(l\Delta t,x)(\phi_i)_x\mathrm{d}x \bigg|
\end{multline}
for $l=1,2,...,100$ where $\widehat{u}_h=\widehat{u}_h(l\Delta t,x;\omega_m)$. This approach could involve training 100 FEONet models over the interval $[0,1]$ by dividing it into 100 segments. However, for efficiency, we divided the interval into 10 segments and modified a model with $u(t=0,x)$ as input and $[u(t=k\Delta t,x)]_{k=1}^{10}$ as output, and then used the output $u(t=10\Delta t,x)$ as input for the next model. Therefore, this approach involved using 10 separate FEONet models.
}

\subsection{Plots on convergence rates}\label{appendix:subsec:convergence}
Let us make some further comments on the experimental details for Figure \ref{fig:rate}. For 1D case, we consider the convection-diffusion equation
\begin{equation}\label{eq_convergence_1d}
  \begin{aligned}
    -\varepsilon u_{xx}+bu_x&=f(x), && x \in D=[-1,1],\\
    u(x)&=0\quad, && x \in \partial D,
  \end{aligned} 
\end{equation} 
where $\varepsilon=0.1$ and $b=-1$. For the 2D case, we consider the following equation
\begin{equation}\label{eq_convergence_2d}
  \begin{aligned}
    -\varepsilon \Delta u+\boldsymbol{v}\cdot \nabla u&=f(x,y), && (x,y) \in D=[-1,1]^2,\\
    u(x,y)&=0, && (x,y) \in \partial D,
  \end{aligned} 
\end{equation}
where $\varepsilon=0.1$ and $\boldsymbol{v}=(-1,0)$. The FEONet was used to approximate the solution of these two equations using P1 and P2 nodal basis functions, and the experiments were conducted with varying domain triangulation to have different numbers of elements to observe the convergence rate. As shown in Figure \ref{fig:rate}, the observed convergence rate of the FEONet shows convergence rates close to 2 and 3 for P1 and P2 approximation respectively, which are the theoretical results for the classical FEM. Since the precision scale of the neural network we used here is about $10^{-2}$, we cannot observe the exact trends in rates below this level of error. We expect to see the same trend even at lower errors if we use larger scale and more advanced models than the CNN structure.

\subsection{Additional plot on generalization capacity}\label{appendix:subsec:generalization}
\begin{figure}[htbp]
\begin{center}
\includegraphics[width=\textwidth]{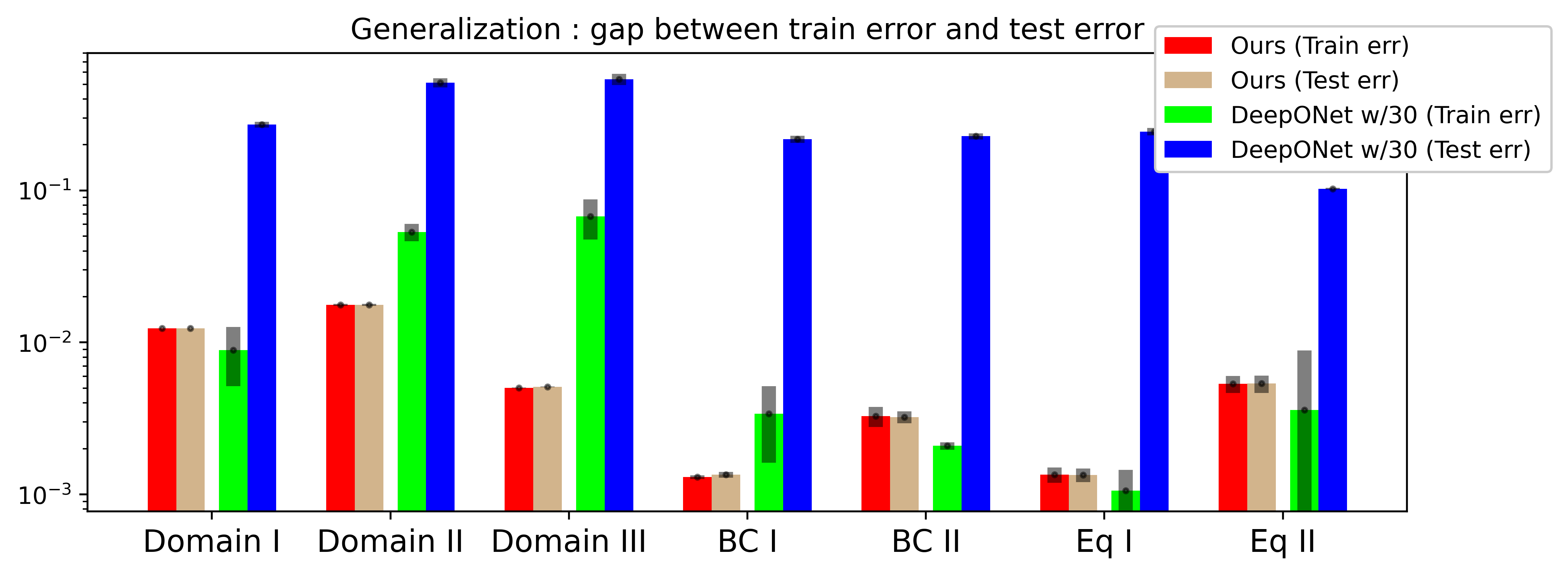}
\end{center}
\caption{Rel. $L^2$ train and test error comparison between the FEONet and the DeepONet with 30 training data. Results are averaged over five independent training trials, and black bars indicate the standard deviation.}
\label{fig:generalization}
\end{figure}
Figure \ref{fig:generalization} displays a comparison of the training and test errors between FEONet (w/o labeled data) and DeepONet (with 30 training data pairs) across different settings as performed in Section \ref{sec:exp} (Table \ref{table:total_error}). We used 3000 input samples to train the FEONet and 30 input-output data pairs to train the DeepONet. Even with the 30 data pairs, there are a lot of computational costs to gain the data. Although it is predictable that increasing the number of input-output data pairs for the DeepONet would narrow the gap between train and test errors, preparing a large amount of data causes significant computational costs. As depicted in Figure \ref{fig:generalization}, the generalization error of FEONet is notably smaller in comparison to that of DeepONet when trained with 30 data pairs.

\section*{Acknowledgments}
Y. Hong was supported by 
Basic Science Research Program through the National Research Foundation of Korea (NRF) funded by the Ministry of Education (NRF-2021R1A2C1093579) and by the Korea government(MSIT) (RS-2023-00219980).
Jae Yong Lee was supported by Institute for Information \& Communications Technology Planning \& Evaluation (IITP) through the Korea government
(MSIT) under Grant No. 2021-0-01341 (Artificial Intelligence Graduate School Program (Chung-Ang University)). S. Ko was supported by National Research Foundation of Korea Grant funded by the Korean Government (RS-2023-00212227).
The authors thank the anonymous referees for their helpful comments that improved the quality of the manuscript.
% We would like to acknowledge the assistance of volunteers in putting
% together this example manuscript and supplement.

\bibliographystyle{abbrv}
\bibliography{references}

\end{document}